\documentclass[a4paper,reqno]{amsart}

\usepackage{natbib}
\usepackage{graphicx}

\usepackage{amssymb}
\usepackage[format=plain,font=small,width=0.3\textwidth]{caption}
\usepackage{amsmath,amsthm}
\usepackage{color}
\usepackage{wrapfig}

\usepackage{ulem} %to use STRIKE OUT command \sout

\makeatletter
\def\@settitle{\begin{center}\baselineskip14\p@\relax\bfseries{\large\@title}\thispagestyle{empty}\end{center}}
\def\@setauthors{%
  \begingroup
  \def\thanks{\protect\thanks@warning}%
  \trivlist
  \centering\footnotesize \@topsep30\p@\relax
  \advance\@topsep by -\baselineskip
  \item\relax
  \author@andify\authors
  \def\\{\protect\linebreak}%
  {\authors}%
  \ifx\@empty\contribs
  \else
    ,\penalty-3 \space \@setcontribs
    \@closetoccontribs
  \fi
  \endtrivlist
  \endgroup
}
\def\maketitle{\par
  \@topnum\z@ % this prevents figures from falling at the top of page 1
  \@setcopyright
  \thispagestyle{firstpage}% this sets first page specifications
  \ifx\@empty\shortauthors \let\shortauthors\shorttitle
  \else \andify\shortauthors
  \fi
  \@maketitle@hook
  \begingroup
  \@maketitle
  \toks@\@xp{\shortauthors}\@temptokena\@xp{\shorttitle}%
  \toks4{\def\\{ \ignorespaces}}% defend against questionable usage
  \edef\@tempa{%
    \@nx\markboth{\the\toks4
      \@nx%\MakeUppercase
      {\the\toks@}}{\the\@temptokena}}%
  \@tempa
  \endgroup
  \c@footnote\z@
  \@cleartopmattertags
}
\makeatother

\newcommand{\out}[1]{}

\newtheorem{Proposition}{Proposition}[section]
\newtheorem{Theorem}[Proposition]{Theorem}

\newtheorem{Corollary}[Proposition]{Corollary}

\theoremstyle{definition}
\newtheorem{Example}[Proposition]{Example}
\newtheorem{Definition}[Proposition]{Definition}
\newtheorem{Remark}[Proposition]{Remark}

\newcommand{\set}[2]{\left\{#1 ;\; #2 \right\}}
\newcommand{\Vin}[1][n]{Q_{G_#1}}
\newcommand{\dVin}[1][n]{P_{#1}}
\newcommand{\pmat}[1]{\begin{pmatrix}#1\end{pmatrix}}
\newcommand{\bs}[1]{\boldsymbol{#1}}
\newcommand{\R}{\mathbb{R}}

\newcommand{\C}{\mathbb{C}}

\newcommand{\ul}[1]{\underline{#1}}

\newcommand{\Disc}{\mathop\mathrm{Disc}}

\newcommand{\supp}{\mathop\mathrm{supp}}
\newcommand{\tlog}[1][a]{\log^{\langle#1\rangle}}
\newcommand{\iunit}{i}

\newcommand{\ds}{\displaystyle}
\newcommand{\pg}[1]{\textcolor{red}{#1}}
\newcommand{\bl}[1]{\textcolor{blue}{#1}}
\newcommand{\la}[1]{\textcolor{magenta}{#1}}

\begin{document}
\title{
%Eigenvalue Distributions of Wigner and Wishart Ensembles of Vinberg Matrices
Wigner and Wishart Ensembles for graphical models
}


\author{Hideto Nakashima}
\address[Hideto Nakashima]{Graduate School of Mathematics, Nagoya University, Furo-cho, Chikusa-ku, Nagoya 464-8602, Japan}
\email{h-nakashima@math.nagoya-u.ac.jp}  
\author{Piotr Graczyk}
\address[Piotr Graczyk]{Laboratoire de Math{\'e}matiques LAREMA, Universit{\'e} d'Angers 2, boulevard Lavoisier, 49045 Angers Cedex 01, France}
\email{piotr.graczyk@univ-angers.fr}

\maketitle

\begin{abstract}
Vinberg cones and the ambient vector spaces are important in modern statistics of sparse models
and of graphical models.
The aim of this paper is to study eigenvalue distributions of 
Gaussian, Wigner and covariance matrices related to
growing Vinberg matrices,
 corresponding  to growing daisy graphs.
For Gaussian or Wigner ensembles,
we give an explicit formula for the limiting distribution.
For Wishart ensembles defined naturally on Vinberg cones,
their limiting Stieltjes transforms, support and atom at 0 are described explicitly 
in terms of the Lambert-Tsallis functions, 
which are defined by using the Tsallis $q$-exponential functions.  
%\keywords{
%Eigenvalue distributions \and
%graphical models
%\and
%covariance matrices \and Wigner matrices 
% \and homogeneous cones \and Vinberg cones\and $q$-exponential  \and Lambert-Tsallis functions
%}
\end{abstract}

\section{Introduction}
\label{sec:intro}

This paper is a first step towards studying {high-dimensional asymptotics of} eigenvalue distributions of 
Gaussian and covariance matrices related to growing  statistical graphical models.

Graphical models provide one of the most powerful methods of unsupervised learning and sparse modelization of modern Data Science and high dimensional statistics 
(cf.\ \citet{Lauritzen,Maat}). 
Mathematical bases of Wishart distributions on matrix cones
related to decomposable and homogeneous graphs 
considered in this paper
were laid down 
by \citet{Lauritzen,LM,Ishi2011LN,GI2014}.

Asymptotics of empirical  eigenvalue distributions are a classical topic of
the random matrix theory (RMT).
There are numerous interactions of RMT with important areas of modern multivariate
statistics: 
high dimensional statistical inference, 
estimation of large covariance matrices, 
principal component analysis (PCA),
time series and many others, 
see the review papers by 
\citet[Section 2]{Diaconis},
\citet{John},
\citet{Paul},
\citet{Bun}, 
the book of \citet{Yao} and the references therein.
%%%%%%%%%%%%%%%%%%%%%%%%%%%%%%%%%%
RMT is also used in signal processing (including MIMO) and compressed sensing
(see \citet[Chapter 10]{Hastie}, for example) 
in the restricted isometry property (RIP) introduced by \citet{CTao}. 
\citet{Fuji} and \citet{Bai} used RMT methods to study  consistency of the  criteria AIC and BIC in estimation of the number of components in PCA.
Distribution of the largest eigenvalue of a Wishart matrix
was studied in \citet{Taka}.

High-dimensional spectral asymptotics for graphical models 
seem to have never been studied before 
and we are convinced that our results will be useful in modern multivariate statistical analysis in the context of graphical models.  
In this paper, 
we concentrate on proving fundamental theorems of RMT,  
the Wigner and Marchenko-Pastur type limit theorems for considered graphical models.
We expect to study statistical applications to estimation of large covariance matrices, the number of
significative PCA factors
and asymptotics of the largest eigenvalue of a sparse Wishart matrix in our subsequent researches.

Growing daisy graphs are among the most natural classes of graphical models. Vinberg matrices are the symmetric matrices corresponding to the growing daisy graphs.
 Covariance matrices are defined naturally on them by a quadratic construction (see Section \ref{ssect:WQ}), thanks to quadratic triangular group actions on 
positive definite Vinberg matrices
(cf.\ Section \ref{sec:Vinberg}). 

In Sections~\ref{sect:gaussian} and~\ref{sect:Wishart},
we provide a complete study of limiting eigenvalue distributions related to Vinberg matrices.
%%%%%%%%%%%%%%%%
The main results are contained in
Theorem \ref{theo:LEDforGW} 
for the Wigner Ensembles
and in
Theorem \ref{th:WishartProfile}
and
Corollaries
\ref{cor:dykh type density},
\ref{prop:support}
and \ref{mainWishart}
for the Wishart Ensembles of Vinberg matrices.
We are able to treat both real and complex matrix ensembles,
but in  view of statistical applications,
we focus on real random matrices.

As a special case of 
Corollary ~\ref{cor:dykh type density},
we provide an elementary and short proof of a result of \citet[{\S}8]{DykemaHaagerup}
on the asymptotic empirical eigenvalue distribution $\mu_0$ 
for the covariance of the triangular real Gaussian  ensemble. 
The proof in~\citet{DykemaHaagerup} is based on the theory of free probability with involved calculations, 
and the Stieltjes transform $S_0(z)$ is given implicitly by determining all the moments of $\mu_0$. 
Later, ~\cite{Cheliotis} mentioned that 
$S_0(z)$ can be 
%written explicitly by using the
expressed in terms of the Lambert  $W$ function.

%%%%%%%%%%%%%%%
Our paper  contributes to the study of triangular random matrices
initiated by \citet{DykemaHaagerup} and continued in \citet{Cheliotis}, 
also in the framework of the theory of Muttalib-Borodin biorthogonal ensembles 
(see \citet{Borodin, Muttalib, Forrester, ForrWang}).
This is a part of recent developments in the theory of singular values of non-symmetric random matrices (see the survey by \citet{Chafai}). 
In contrast to \citet{Cheliotis},
we do not dispose of an explicit formula for the joint eigenvalue density.

The analysis, 
probability and statistics on homogeneous cones develops intensely in recent years (\citet{AW,GI2014,GIL2019,Ishi2011LN,Ishi2016,LM,YN2015,N2020}),
and
Vinberg cones and dual Vinberg cones are basic examples of homogeneous cones
(see Section \ref{sec:Vinberg}).
Our results are   a first contribution to the RMT on homogeneous cones. 

The main method used in our paper is the {\it variance profile method} 
for Gaussian and Wigner matrix ensembles, 
presented in Section~\ref{sect:2-5}.
It was
applied first in \citet{Shly} in the Gaussian case
and developed in \citet{AndersonZeitouni2006}
in the Wigner case. 
We use
the recent approach of \citet{LNofB}.
In Theorem~\ref{theo:variance profile method}
we slightly strengthen for our needs the main variance profile result of \citet{LNofB}.
Theorem~\ref{theo:variance profile method}
will be  useful for studying of eigenvalue distributions 
related to  general graphical models. 

Note that the variance profile methods were also developed directly  for Wishart ensembles by
\citet{HLN,HLN-AHP,HLN-AAP214,HLN-AAP515} 
(cf.\ Remark~\ref{otherWishart}).
The variance profile methods are related to operator-valued free probability theory 
(\citet[Chapter 9]{MingoSpeicher}).

Our expression of a limiting
Stieltjes transform 
for Wishart Ensembles of Vinberg matrices,
is based on the introduction of Lambert-Tsallis functions $W_{\kappa,\gamma}$,
see Section~\ref{sect:functionW}.
The Lambert-Tsallis functions are defined by using  Tsallis $q$-exponential functions,
now actively studied in Information Geometry (cf.\ \citet{AmariOhara2011,ZNS2018}).

Outlines of all proofs are given. Technical details 
are omitted and  can be viewed
in Supplementary material available
from the editor of the journal.

Simulations of histograms of eigenvalues of Vinberg matrices
are illustrated by
Figures \ref{fig:1}-\ref{fig:5} 
in the Wigner case and by Figures \ref{Wishart2}-\ref{Wishart4}
in the Wishart case.

\section{Preliminaries}
\label{sect:preliminaries}

We begin this paper with recalling the definition of the empirical eigenvalue distribution of      a symmetric matrix.
Let $X\in\mathrm{Sym}(n,\R)$ be a symmetric matrix
and let {$\lambda_1(X)\ge \cdots\ge \lambda_n(X)$ be the ordered} eigenvalues of $X$ with counting multiplicities. Denote by $\delta_a$  the Dirac measure at $a$.
Then, the empirical eigenvalue distribution $\mu_X$ of $X$ is defined by
$
\mu_X=\frac{1}{n}\sum_{i=1}^n\delta_{\lambda_i(X)}.
$

If $\{X_n\}_{n=1}^{\infty}$ $(X_n\in\mathrm{Sym}(n;\,\R))$ is a sequence of Gaussian, Wigner or Wishart matrices,
then it is well known that there exists a limit $\mu$ of $\mu_{X_n}$ as $n\to\infty$, and
the sequence of random measures $\mu_{X_n}$ converges almost surely weakly  to
 the semi-circle law or the Marchenko-Pastur law, respectively
(see for example \cite{BaiSilverstein,LNofB}).
The limits $\mu$ of $\mu_{X_n}$,
in the almost sure weak sense, are  said to be the ``limiting eigenvalue distributions $\mu$ of $X_n$.''  For simplicity,
we will say ``i.i.d.\ matrices'' instead of  ``matrices with independent and identically distributed  non-null terms''.

%In this paper,
%we consider random matrices on which we impose some sparsity restrictions
%{ that arise} naturally from the  graphical models.

%%%%%%%%%%%%%%%%%%%%%%%%%%%%%%%%%%%%%%%%%%%%%%%%%%%%%%%%%%%%%%%%%%%%%%%%%%%%%%%%%%
\subsection{Basics on statistical graphical models}
\label{ssect:graph}

Let $G$ be a graph with vertices $V=\{1,2,\ldots, n\}$ and edges $E$.
We say that a statistical character $\mathcal X=(X_1,\ldots,X_n)$
has the {\it dependence graph} $G$ when
each conditional independence of marginals $X_i$ and $X_j$ with respect to remaining variables corresponds to the absence of the edge $\{i,j\}$ in $E$. 
Thus the dependence graph $G$ is a  tool of encoding of the
conditional independence of marginals of ${\mathcal X}$.
We say that ${\mathcal X}$ belongs to the graphical model governed by $G$.

Let $\mathcal{U}_G$ be the subspace of $\mathrm{Sym}(n,\R)$ containing matrices with $u_{ij}=0$
if the edge $\{i,j\}\not\in E$.
Cones $P_G=\mathrm{Sym}(n,\R)^+\cap \mathcal{U}_G$
and their dual cones $Q_G$ are basic objects of graphical model theory. 
Actually, 
a Gaussian $n$-dimensional model $N(m,\Sigma)$ is governed by the graph $G$
if and only if the inverse covariance matrix $\Sigma^{-1}\in P_G$
(cf.\ \cite{Lauritzen}).

\begin{wrapfigure}[11]{r}{100pt}
\vspace{-2em}
\centering
\includegraphics[width=100pt,height=100pt]{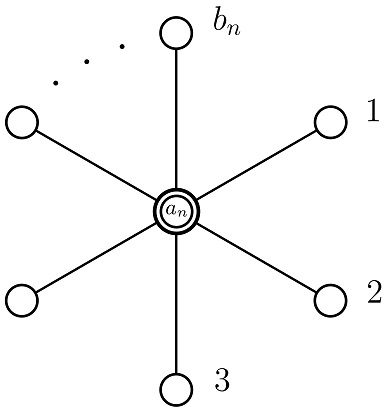}
\caption{Daisy Graph}
\label{fig:Daisy Graph}
\end{wrapfigure}

An important class of graphical models, called {\it daisy graphs}, is defined as follows.
Let $a+b=n$ and let $D(a,b)$ be a graph with vertices $V =\{1,\dots,n\}$,
such that the first $a$ elements form a complete graph and 
the latter $b$ elements are satellites (petals) of the complete graph, 
that is, each satellite connects to all elements in the complete graph 
and does not connect to the other satellites (see Figure~\ref{fig:Daisy Graph}).
{The  double circle  around the vertex  $a_n$ in Figure~\ref{fig:Daisy Graph} indicates the complete graph with $a_n$ vertices.}

In high dimensional statistics, 
it is essential to let the number of observed characters $n$ tend to infinity. 
From the graphical model theory point of view,  
the pattern of the growing graphs $G_n$ and of the corresponding cones $P_{G_n}$ 
should remain the same.
This requirement is met by growing daisy graphs $D(a_n,b_n)$ for non-decreasing sequences of positive integers
$\{a_n\}_{n=1}^{\infty}$ and $\{b_n\}_{n=1}^{\infty}$ such that $a_n+b_n=n$. 

\subsection{Generalized dual Vinberg cones and Vinberg matrices} 
\label{sec:Vinberg}

Let $\{a_n\}_{n=1}^{\infty}$ and $\{b_n\}_{n=1}^{\infty}$ be non-decreasing sequences of positive integers 
such that $a_n+b_n=n$ and the ratio $a_n/n$ converges to $c\in[0,1]$.
Let $G_n=D(a_n,b_n)$ be the corresponding daisy graph. 
Then, the corresponding matrix space $\mathcal{U}_n$ of the graph $G_n$ is a subspace of $\mathrm{Sym}(n,\R)$ defined by
\[
\mathcal{U}_n:=\set{{U}=\pmat{x&y\\ {}^{t}y&d}}{\begin{array}{l}x\in\mathrm{Sym}(a_n,\R),\,y\in\mathrm{Mat}(a_n\times b_n,\R),\\\text{$d$ is a diagonal matrix of size $b_n$}\end{array}},
\]
and we set
%\begin{equation} \label{def:dual cone}
\[
\dVin:=P_{G_n}=\mathcal{U}_n\cap \mathrm{Sym}(n,\R)^+.
\]
%\end{equation}
Then, $\dVin$ is an open convex cone in $P_{G_n}$.
Moreover, the cone $\dVin$ admits a transitive group action, 
{\it i.e.} $\dVin$ is a {\it  homogeneous cone},
since the following triangular group
\[
H_n:=\set{h=\pmat{h_1&y\\0&d}\in GL(n,\R)}{\begin{array}{l}
h_1\in GL(a_n,\R)\text{ is upper triangular},\\
y\in\mathrm{Mat}(a_n\times b_n;\,\R),\\ d\colon\text{diagonal of size }b_n
\end{array}}
\]
acts on $\dVin$ transitively by the
quadratic action $\rho(h){U}:=h{U}{}^{\,t\!}h$ for $h\in H_n$ and ${U}\in\dVin$.
This is easily verified by using the Cholesky decomposition
(cf.\ \citet[p.\ 3]{Ishi2016}).
For definition and basic properties of homogeneous cones, see \citet{Vinberg, Ishi2011LN}.

If $n=3$ and $(a_n,b_n)=(1,2)$, 
then $\dVin[3]$ is the dual Vinberg cone
(see Example \ref{ex:Vinberg})
so that, in this paper, 
we call $\dVin$ a \textit{generalized dual Vinberg cone}
and elements $U\in\mathcal{U}_n$ \textit{Vinberg matrices}.
Vinberg cones form an important class
of matrix cones related to  graphical models (cf.\ Section \ref{ssect:graph}).
On the other hand,
if we set $a_n=n-1$ and $b_n=1$, 
then $\mathcal{U}_n$ is the space $\mathrm{Sym}(n,\R)$ of symmetric matrices of size $n$,
and hence our discussion covers the classical results.
In what follows,
we introduce two kinds of random matrices related to the homogeneous cones $P_n$, that is, Gaussian and  Wigner matrices   and Wishart quadratic (covariance) 
matrices.
%%%%%%%%%%%%%%%%%%%%%%%%%%%%%%%
\subsection{Gaussian  and Wigner matrices in  $\mathcal{U}_n$}
\label{ssect:GW}
%%%%%%%%%%%%%%%%%%%%%%%%%%%%%
Analogously to the classical Wigner matrices,
we say that $U_n=(u_{ij})\in\mathcal{U}_n$ is
a Wigner  random matrix if
\begin{equation}
\label{def:gaussianU}
    \text{
    $\left\{\rule{0pt}{40pt}\right.$
    \begin{tabular}{@{}l@{\ }p{0.8\textwidth}}
    $\bullet$&the diagonal terms $(u_{ii})$ are independent of the off-diagonal terms $(u_{ij})_{i<j}$, \\
    $\bullet$&the diagonal $u_{ii}$'s are  centered i.i.d. variables with variance $v'$ and
    fourth moment $M_4'$,\\
    $\bullet$&the non-nul off-diagonal $u_{ij}$'s, $i<j$, are   centered i.i.d. variables with variance $v$
    and fourth moment ${M_4}$,
    \end{tabular}}
\end{equation}
where $v,v',M_4,M_4'$ are fixed positive real numbers.
If the non-nul terms $u_{ij}$
are Gaussian, with $\nu=1$ and $\nu'=2$, the matrices $U_n$ form a Gaussian Orthogonal Ensemble of Vinberg matrices.

In Section \ref{sect:gaussian}, we consider empirical eigenvalue distributions of rescaled Wigner matrices $U_n/\sqrt{n}\in\mathcal{U}_n$.
%%%%%%%%%%%%%%%%%%%%%%%%%%%%%%%%%%
\subsection{Quadratic construction of   Wishart  (covariance) matrices in  $\mathcal{U}_n$}
\label{ssect:WQ}
Recall that Wishart matrices are constructed quadratically both in Random Matrix Theory and in statistics.
In this section we define, by a
quadratic construction,
Wishart  (covariance) matrices in $\mathcal{U}_n$.

We first recall the notion of a direct sum of quadratic maps.
Let $Q_i\colon\R^{m_i}\to\R^{m}$ $(i=1,\dots,k)$ be quadratic maps.
Then, the direct sum $Q_1\oplus\cdots\oplus Q_k$ is an $\R^{m}$-valued quadratic map on $\R^{m_1}\oplus\cdots\oplus\R^{m_k}$ given by
\[
Q(x):=Q_1(x_1)+\cdots+Q_k(x_k)\quad\text{where}\quad
x=\sum_{i=1}^kx_i\quad\bigl(x_i\in\R^{m_i}\bigr).
\]
If $Q_1=\cdots=Q_k$, then the direct sum $Q$ is denoted by $Q_1^{\oplus k}$.
As showed in \citet{GI2014},
any homogeneous cone $\Omega$ admits 
a canonical  family of the so-called \textit{basic quadratic maps} 
$q_j$ ($j=1,\ldots, r$)
defined for each $j$ on a suitable finite dimensional vector space $E_j$
and with values in the closure $\overline{\Omega}$ of $\Omega$.
The number $r$ is called the rank of $\Omega$ 
and $r=n$ for
the cones $\mathcal{U}_n$.
Using the basic quadratic maps $q_j$, one constructs  quadratic maps $Q_{\ul{k}}$ for $\ul{k}\in\mathbb{Z}^{r}_{\ge0}$ by
\[Q_{\ul{k}}:=q_1^{\oplus k_1}\oplus\cdots\oplus q_r^{\oplus k_r},\]
defined on $E_{\ul{k}}:=E_1^{\oplus k_1}\oplus\cdots\oplus E_r^{\oplus k_r}$.
The maps $Q_{\ul{k}}$ are $\Omega$-positive, {\it i.e.}
if $\xi\in E_{\ul{k}}\setminus\{\bs{0}\}$, then
$Q_{\ul{k}}(\xi)\in\overline{\Omega}\setminus\{\bs{0}\}$.

In our case $\Omega=\dVin$, the basic quadratic maps are given as follows
(cf.\ \cite{GI2014}).
For $j=1,\dots,n$, define
$E_j\subset \R^n$ by
\[
\begin{array}{l}
E_j=\set{\pmat{\bs{\xi}\\ \bs{0}}\in\R^n}{\bs{\xi}\in\R^j}\quad(j\le a_n),\\[0.8em]
E_j=\set{\pmat{\bs{\xi}\\ \bs{0}}+\xi'\bs{e}_{j}\in\R^n}{\bs{\xi}\in\R^{a_n},\ \xi'\in\R}\quad(j>a_n),
\end{array}
\]
where $\bs{e}_{i}$ $(i=1,\dots,n)$ is the  vector in $\R^{n}$ having $1$ on the $i$-th position and zeros elsewhere.
We note that each $E_j$ corresponds to the $j$-th column of the Lie algebra $\mathfrak{h}_n$ of $H_n$,
that is,
we have
$\mathfrak{h}_n=\set{H=(\bs{\xi}_1,\dots,\bs{\xi}_n)}{\bs{\xi}_j\in E_j}$.
Then,
the basic quadratic maps $q_j\colon E_j\to \mathcal{U}_n$ of the cone $\dVin$ are defined by
\[
q_j(\bs{\xi}_j):=\bs{\xi}_j{}^{\,t\!}\bs{\xi}_j\in\mathcal{U}_n\quad(\bs{\xi}_j\in E_j).
\]
Let $\ul{k}\in\mathbb{Z}_{\ge0}^n$.
Then, $E_{\ul{k}}$ can be viewed as a subspace of $\mathrm{Mat}(n\times(k_1+\cdots+k_n);\,\R)$.
In fact, we have
\[
\begin{array}{l}
\ds
E_{\ul{k}}=\set{\eta=\Bigl(\overbrace{\bs{\xi}_1^{(1)},\dots,\bs{\xi}_1^{(k_1)}}^{k_1},\bs{\xi}_2^{(1)},\dots,\bs{\xi}_{n-1}^{(k_{n-1})},\overbrace{\bs{\xi}_n^{(1)},\dots,\bs{\xi}_n^{(k_n)}}^{k_n}\Bigr)}{\begin{array}{l}
\bs{\xi}_j^{(i)}\in E_j,\\
j=1,\dots,n,\\
i=1,\dots,k_j
\end{array}
}\\[1em]
\hfill{}
\subset \mathrm{Mat}(n\times(k_1+\cdots+k_n);\,\R),
\end{array}
\]
and then $Q_{\ul{k}}(\eta)=\eta{}^{\,t\!}\eta$ for $\eta\in E_{\ul{k}}$.

When $\eta\in E_{\ul{k}}$ is an i.i.d.\ random matrix whose non-null
terms have the normal law $N(0,v)$,
the law of  $Q_{\ul{k}}(\eta)$ 
is a
{\it Wishart law} $\gamma_{Q_{\ul{k}},1/(2v){\rm Id}_n}$
on the cone $P_n$. 
 For the definition of all 
 Wishart laws on the cone $P_n$, see \cite{GI2014}.
More generally, in this paper,
we consider eigenvalue distributions of rescaled matrix $Q_{\ul{k}}(\eta)/n$ 
under the assumption that
$\eta\in E_{\ul{k}}$ is a centered rectangular i.i.d.\ matrix whose non-null
terms have variance $v$ 
and finite fourth moments $M_4$.

We   consider two-dimensional  multiparameters $\ul{k}=\ul{k}(n)\in\mathbb{Z}_{\ge0}^n$  of the form
\begin{equation}\label{def:k}
\ul{k}=m_1(1,\dots,1)+m_2(\,\overbrace{0,\dots,0}^{a_n},\overbrace{1,\dots,1}^{b_n}\,)\quad(m_1,m_2\in\mathbb{Z}_{\ge0}).
\end{equation}

\begin{Example}\label{ex:Vinberg}
{\rm
Let $n=3$, $a_3=1$ and
$b_3=2$.
In this case,
$\dVin[3]$ is the  dual Vinberg cone 
(cf.\ \citet[p.\ 397]{Vinberg}, \citet[ {\S}5.2]{Ishi2001}):
\[
\dVin[3]=\set{x=\pmat{x_{11}&x_{12}&x_{13}\\x_{12}&x_{22}&0\\x_{13}&0&x_{33}}}{\textrm{$x$ is positive definite}}.
\]
Consider $m_1=m_2=1$, so  $\ul{k}=(1,2,2)$. Then $E_{\ul{k}}=E_{(1,2,2)}$ can be written as
\[
E_{(1,2,2)}=\set{
\eta=
\pmat{x&y_{11}&y_{12}&z_{11}&z_{12}\\0&y_{21}&y_{22}&0&0\\0&0&0&z_{21}&z_{22}}}{x,y_{ij},z_{ij}\in\R},
\]
and $Q_{(1,2,2)}(\eta)=\eta{}^{\,t\!}\eta$ is given as
\[
Q_{(1,2,2)}(\eta)=
\pmat{x^2+y_{11}^2+y_{12}^2+z_{11}^2+z_{12}^2&y_{11}y_{21}+y_{12}y_{22}&z_{11}z_{21}+z_{12}z_{22}\\
y_{11}y_{21}+y_{12}y_{22}&y_{21}^2+y_{22}^2&0\\
z_{11}z_{21}+z_{12}z_{22}&0&z_{21}^2+z_{22}^2
}.
\]
If $x,y_{ij},z_{ij}$ are $N(0,v)$ i.i.d. Gaussian variables, the random matrix $Q_{(1,2,2)}(\eta)$ has a Wishart law on $P_3$.
}
\end{Example}

The form \eqref{def:k} of the Wishart multiparameter $\ul{k}$
englobes and generalizes the following cases.
In both cases, with rescaling 
$1/n$, the limiting  eigenvalue distribution is known. 

\begin{enumerate}
    \item[(i)]
    The classical Wishart Ensemble $M\,^t M $ on $\mathrm{Sym}(n,\R)^+$,
    where $M=M_{n\times N}$ is an i.i.d.\ matrix with finite fourth moment ${M_4}$, 
    with parameter $C:=\lim_n\frac{N}{n}>0$
    (see \citet{AGZ,Faraut})  
    for $(a_n,b_n)=(n-1,1)$, $m_1=0$ and $m_2\sim Cn$.
    The limiting eigenvalue distribution is the Marchenko-Pastur law $\mu_C$ with parameter $C$,
    i.e.\ denoting $a=\bigl(\sqrt{C}-1\bigr)^2,
    b=\bigl(\sqrt{C}+1\bigr)^2$
    and $[x]_+:=\max(x,0)$ $(x\in\R)$,
    \[
    \mu_C=[1-C]_+\delta_0+\frac{\sqrt{(t-a)(b-t)}}{2\pi t}\,\chi_{[a,b]}(t)dt.
    \]
    \item[(ii)]
    The  Wishart Ensemble related to the Triangular Gaussian Ensemble\\ 
    (\citet{DykemaHaagerup, Cheliotis}) 
    for $(a_n,b_n)=(n-1,1)$, $m_1=1$  and $m_2=0$.
    When $v=1$,
    the limiting eigenvalue distribution, 
    which we call the {\it Dykema-Haagerup measure} $\chi_1$, 
    is absolutely continuous with respect to Lebesgue measure and has support equal to the interval $[0,e]$.
    Its density function $\phi$ is defined on the interval $(0,e]$ by the implicit formula
    (\citet[Theorem 8.9]{DykemaHaagerup})
    \begin{equation}\label{dykh}
    \phi\left(\frac{\sin x}{x}\,\exp(x\cot x)\right)
    =
    \frac{1}{\pi}\sin x\exp(-x\cot x)\qquad(0\le x<\pi),
   \end{equation}
    with $\phi(0+)=\infty$
    and $\phi(e)=0$.
    For $v\not=1$, the limiting measure $\chi_v$ has density $\phi({y}/{v})/v$ on
   the segment $(0,ve]$.
\end{enumerate}
\subsection{Resolvent method for 
Wigner ensembles with a variance profile $\sigma$}
\label{sect:2-5}

Let $\C^+$ denote the upper half plane in $\C$.
In this paper, the Stieltjes transform $S(z)=S_{\mu}(z)$ of a probability measure $\mu$ on $\R$ is defined to be
\[
S(z)=\int_{\R}\frac{\mu(dt)}{t-z}\quad(z\in\C^+).
\]

In the sequel, we will need the following properties of the Stieltjes transform, which are not difficult to 
prove.
\begin{Proposition}
 \label{th:0 ONE}
 1. Suppose that $s(z)$ is the Stieltjes transform of a finite  measure $\nu$ on $\R$.
If for all $x\in\R$ it holds
\[
\lim_{y\to 0+} \mathrm{Im}\, s(x+iy)=0
\]
then $ s(z)\equiv 0$ and $\nu$ is a null  measure ($\nu(B)=0$ for any Borel set $B$).
\\
2.  
 Suppose $f\ge 0$ and  $f\in L^1(\R)$.
Let $s(z)$ be the Stieltjes transform of $f$.
If $f$ is continuous at $x$ then
\begin{equation} \label{eq:StContfONE}
\lim_{y\to 0+} \frac1{\pi} \mathrm{Im}\,s(x+iy)=f(x).
\end{equation}
If $f$ is continuous on an interval
$[a,b]$, $a<b$, the convergence 
\eqref{eq:StContfONE}
is  uniform for $x\in[a,b]$.
\end{Proposition}
Recall that if $\mu$ is a probabilistic measure on $\R$, with Stieltjes transform $s(z)$
and the absolutely continuous part of $\mu$ has density   $f$, then
\eqref{eq:StContfONE} holds
 for almost all $x$ (Lemma 3.2 (iii) of \citet{LNofB}).

We
present now the following, slightly 
strengthened result
from the Lecture Notes of \citet[{\S}3.2]{LNofB}, 
 that will be a main tool of proofs in this paper.

Let $\sigma\colon[0,1]\times[0,1]\to [0,\infty)$ be a 
bounded Borel measurable symmetric function.
For each integer $n$,
we partition the interval $[0,1]$ into $n$ 
equal intervals $J_i,
i=1,\ldots, n$.
Put $Q_{ij}:=J_i\times J_j$, which is a partition of $[0,1]\times[0,1]$.
We assume that $Y_{ij}$ $(i\le j)$ are independent centered real variables, defined on a common probability space, with variance
\begin{equation}
\label{eq:varianceprofile}
\mathbb{E}Y_{ij}^2
=
\frac{1}{n}
\left(
\int_{Q_{ij}}
\frac{\sigma(x,y)}{|Q_{ij}|}\,dx\,dy
+\delta_{ij}(n)
\right), 
\end{equation}
for a sequence $\delta_{ij}(n)$.
We note that the law of $Y_{ij}$ depends on $n$.
We set $Y_{ji}:=Y_{ij}$ and we consider the symmetric matrix
$
Y_n:=(Y_{ij})_{1\le i,j\le n}.
$
We note that,
if $\sigma$ is continuous, then, up to a perturbation $\delta_{ij}(n)$, the variance of $\sqrt{n}Y_{ij}$ is 
approximatively ${\sigma(i/n,j/n)}$,
and hence we call {$\sigma$} a variance profile 
in this paper.

\begin{Theorem}
\label{theo:variance profile method}
%{sect:2-5}
Let  $\delta_0(n) := \displaystyle \frac{1}{n^2} \sum_{i,j\le n} |\delta_{ij}(n)|$. 
Assume \eqref{eq:varianceprofile} and suppose that
\begin{equation}
\label{eq:moment}
 \lim_n \delta_0(n)= 0
\quad
\text{and}\quad
 \max_{i,j\le n} 
\frac{\mathbb{E}(Y_{ij}^4)}{n(\mathbb{E}Y_{ij}^2)^2} = o(1)
\quad (Y_{ij}\not=0).
\end{equation}
Let $\mu_{Y_n}$ be the empirical eigenvalue distribution of $Y_n$.
Then, there exists a probability measure $\mu_\sigma$ depending on $\sigma$
such that $\mu_{Y_n}$ converges weakly to $\mu_{\sigma}$ almost surely.
The Stieltjes transform $S_\sigma$ of $\mu_\sigma$ is given as follows.\\
(a) For each $z$ with $\mathrm{Im}\,z>1$,
there exists a
unique $\C^+$-valued $L^1$-solution
 $\eta_z:[0,1]\mapsto \C^+$,
 of the equation
\begin{equation}
\label{eq:FE of vp}
\eta_z(x)=-\left(z+\int_0^1\sigma(x,y)\,\eta_z(y)\,dy\right)^{-1}
\quad(\text{for almost all $x\in[0,1]$}),
\end{equation}
and the function $z\mapsto\eta_z(x)$ extends to an analytic 
$\C^+$-valued function
on $\C^+$, for almost all $x\in[0,1]$.
Then,
\[
S_\sigma(z)=\int_{0}^1\eta_z(x)\,dx.
\]
(b) The function $x\to \eta_z(x) $
is also a solution of  \eqref{eq:FE of vp} for $0<\mathrm{Im}\,z\le 1.$
\end{Theorem}
\begin{proof}
The proof is the same as the proof of \citet[Theorem 3.1]{LNofB}, 
where a stronger assumption $|\delta_{ij}(n)|\le \delta(n)$ is required for some sequence $\delta(n)$ going to 0. 
It is replaced by the first condition of \eqref{eq:moment}.
Detailed   analysis of the  proof of
the approximate fixed point equation  in \citet[page 42]{LNofB}
shows that 
the weakest assumption on the fourth moments $\mathbb{E}Y_{ij}^4$
ensuring the concentration of the conditional variance  of $\langle
Z_i, R^{(i)}Z_i\rangle$ is
the second  condition of \eqref{eq:moment}.
The property (b) is observed in 
\citet[page 39]{LNofB} by analiticity.
\qed
\end{proof}

Theorem~\ref{theo:variance profile method} shows that, 
to each variance profile function $\sigma$, 
one associates uniquely a Stieltjes transform $S_\sigma(z)$
of a probability measure. 
For the correspondence between $\sigma$ and $S_\sigma$, 
the conditions (7) are not needed.
We define $S_\sigma(z)$ as the {\it Stieltjes transform associated to $\sigma$}. 

\begin{Remark}
{\rm
A prototype of the
variance profile method for
Wigner ensembles was given
by~\citet[Theorem 3.2]{AndersonZeitouni2006}.
Theorem 3.1 of \citet{LNofB}
and Theorem \ref{theo:variance profile method}
provide a simple general approach. 
Special cases of
variance profile convergence results
for Wigner matrices
were studied before, as discussed  below in (i) and (ii).

\noindent
\begin{tabular}{c@{\hspace{1ex}}p{0.95\textwidth}}
(i)&
If we set $\sigma(x,y)=1$ for all $x,y$, 
then $\sqrt{n}Y$ is a Wigner ensemble with $v=v'=1$.
Let $S_{\mathrm{sc}}(z)$ be the Stieltjes transform of the semi-circle law on $[-2,2]$.
Then, 
the functions $x\to \eta_z(x)$  do not depend on $x$ (but do on $z$) and 
the functional equation~\eqref{eq:FE of vp} gives 
the equation $S_{\mathrm{sc}}(z)=-(z+S_{\mathrm{sc}}(z))^{-1}$,
which
is well known from the detailed study of resolvent matrices
(see \citet[{\S2}.4.3]{Tao}).\\
(ii)& 
The paper \citet{AndersonZeitouni2006}
deals primarily with a variance profile $\sigma$ such that
$\int\sigma(x,y)\,dy=1$ for any $x$,
corresponding to a band matrix model.
For band matrix ensembles, see also~\citet{ErdosAdM,Erdos,NicaShlyyakhtenkoSpeicher2002,Shly}.
\end{tabular}
}
\end{Remark}

\section{Wigner Ensembles of Vinberg Matrices}
\label{sect:gaussian}

\begin{figure}[ht]
    \centering
    \begin{tabular}{ccc}
    \begin{minipage}{0.33\textwidth}
        \centering
        \includegraphics[scale=0.2]{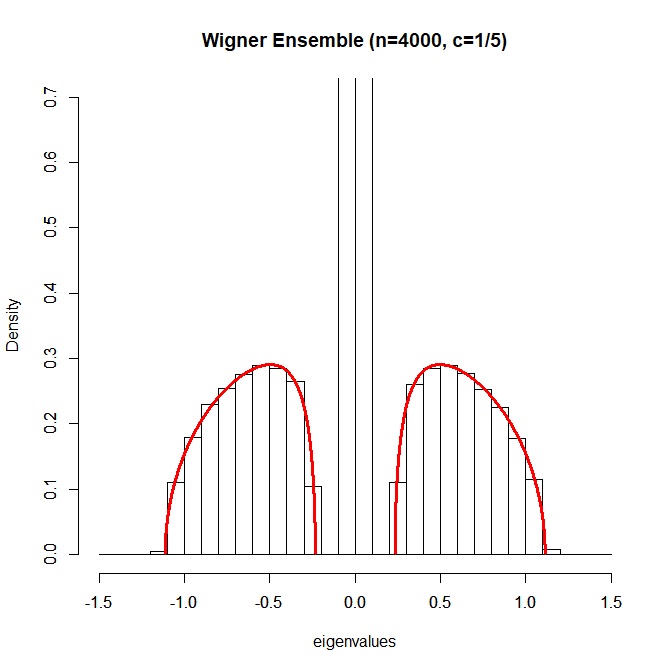}
        \caption[1h]{Simulation for $c=1/5$}
        \label{fig:1}
    \end{minipage}
    &
    \begin{minipage}{0.31\textwidth}
        \centering
        \includegraphics[scale=0.2]{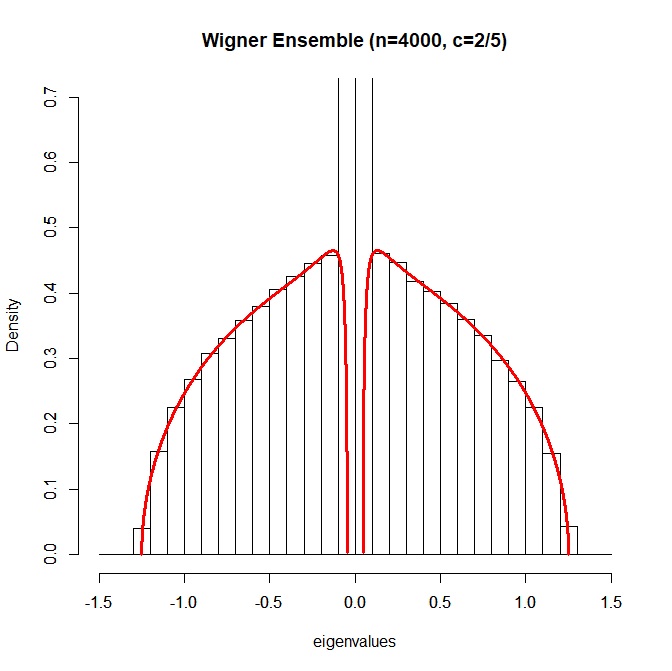}
        \caption{Simulation for $c=2/5$}
        \label{fig:2}
    \end{minipage}
    &
    \begin{minipage}{0.31\textwidth}
        \centering
        \includegraphics[scale=0.2]{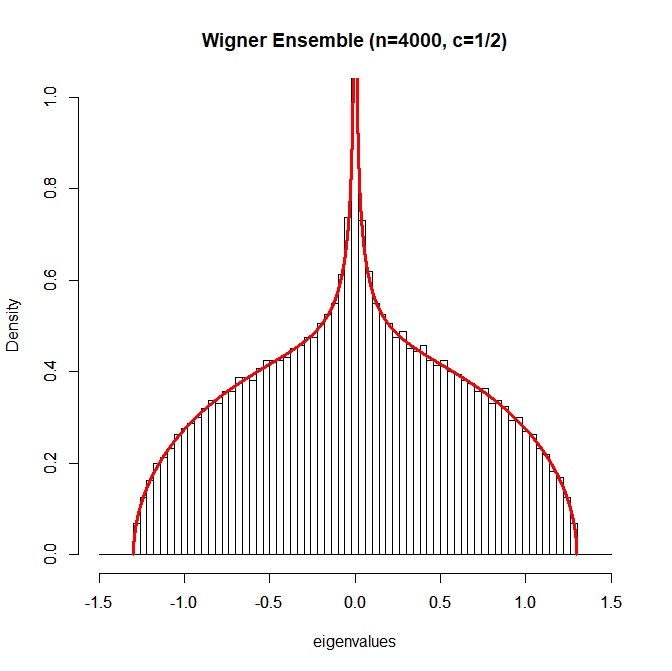}
        \caption{Simulation for $c=1/2$}
        \label{fig:3}
    \end{minipage}\\
    \begin{minipage}{0.31\textwidth}
        \centering
        \includegraphics[scale=0.2]{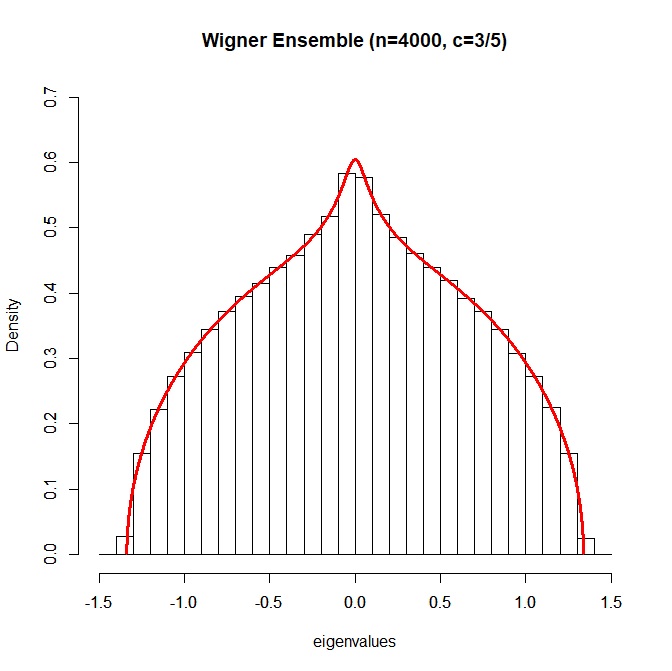}
        \caption{Simulation for $c=3/5$}
        \label{fig:4}
    \end{minipage}
    &
    \begin{minipage}{0.31\textwidth}
        \centering
        \includegraphics[scale=0.2]{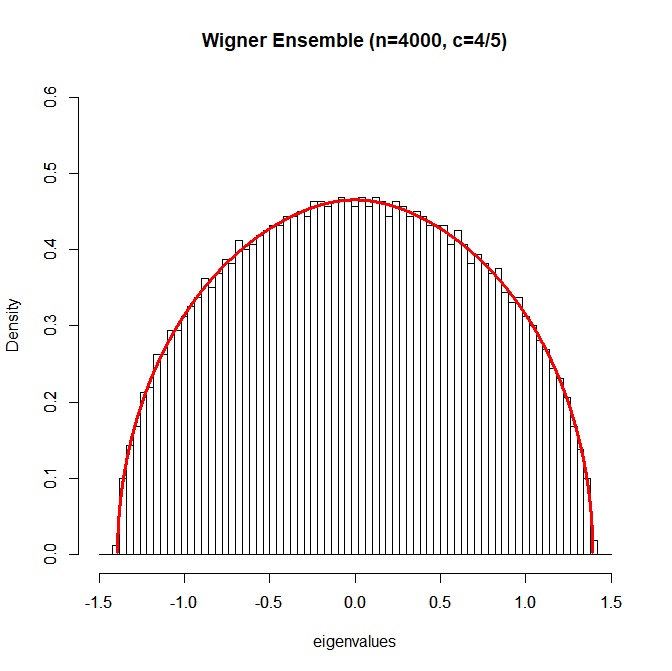}
        \caption{Simulation for $c=4/5$}
        \label{fig:5}
    \end{minipage}
    &
    \begin{minipage}{0.31\textwidth}
        \centering\rule{0pt}{20pt}
        \includegraphics[scale=0.26]{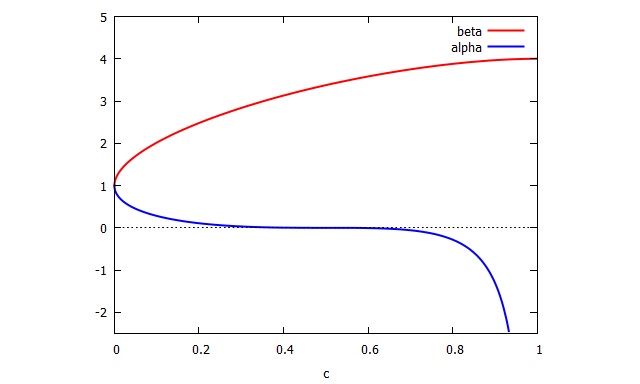}\\[0.5em]
        \caption{Graphs of $\alpha_c$ and $\beta_c$}
        \label{fig:6}
    \end{minipage}
    \end{tabular}
\end{figure}

In this section,
we give explicitly the limiting eigenvalue distributions $\mu$ for the
scaled Wigner matrices $U_n\in\mathcal{U}_n$ defined by \eqref{def:gaussianU}.
Let $\chi_I$ denote the indicator function of a subset $I\subset \R$.
For a real number $a$, its cubic root is denoted by $\sqrt[3]{a}\in\R$
and set $[\,a\,]_+=\max(a,0)$.
We introduce two real  numbers $\alpha_c$, $\beta_c$ depending on $c\in[0,1)$ by
\begin{equation}
\label{eq:alphabeta}
\alpha_c=\frac{8+4c-13c^2-\sqrt{c(8-7c)^3}}{8(1-c)},\ 
\beta_c=\frac{8+4c-13c^2+\sqrt{c(8-7c)^3}}{8(1-c)}.
\end{equation}
It is clear that 
$\alpha_0=\beta_0=1$,  $\alpha_c < \beta_c$ and $\beta_c>0$ for all $c\in(0,1)$.
We note that 
$\alpha_{1/2}=0$, $\alpha_c<0$ when $c>1/2$,
 $\lim_{c\to 1-}\alpha_c=-\infty$,
 $\lim_{c\to1-}(1-c)\alpha_c=-1/4$  and
$\lim_{c\to1-}\beta_c=4$,
 so that we set $ \beta_1=4$. It can be shown that $c\mapsto \alpha_c$ is strictly decreasing and $c\mapsto \beta_c$ is strictly increasing on $[0,1]$
 {(see Figure~\ref{fig:6}).}

\begin{Theorem}
\label{theo:LEDforGW}
Let $U_n$ be a Wigner matrix 
on $\mathcal{U}_n$ defined by \eqref{def:gaussianU}. 
Assume that $\lim_{n\to+\infty}a_n/n=c\in(0, 1)$.
Then, the limiting eigenvalue distribution $\mu$ of the rescaled  matrices $U_n/\sqrt{n}$ exists
and is given for $c\in (0,1)$ as
\[
\mu=f_c(t)\,dt+[1-2c]_+\delta_0
\]
with
\begin{equation}
\label{eq:mu}
f_c(t):=  \frac{\sqrt[3]{R_{+}\left(t/\!\sqrt{v};\,c\right)}
-
\sqrt[3]{R_{-}\left(t/\!\sqrt{v};\,c\right)}}{2\sqrt{3}\pi\,t}\,\chi_{[\alpha_c,\beta_c]}\left(\frac{t^2}{v}\right),
\end{equation}
where,  for $x^2\in [\alpha_c,\beta_c]$,
\[
\begin{array}{r@{\ }c@{\ }l}
R_{\pm}(x;\,c)&
:=&
x^6-3(c+1)x^4+\frac{3}{2}(5c^2-2c+2)x^2+(2c-1)^3\\[0.6em]
& &\qquad\pm
3c\sqrt{3-3c}\cdot x\sqrt{(x^2-\alpha_c)(\beta_c-x^2)}.
\end{array}
\]
The support of $\mu$ is given as
\out{$\{0\}$ when $c=0$,} 
\begin{equation}
\label{eq:supp}
\!\supp \mu=\begin{cases}
\left[\,-\sqrt{v\beta_c},\,-\sqrt{v\alpha_c}\,\right]\cup \{0\} \cup
\left[\,\sqrt{v\alpha_c},\,\sqrt{v\beta_c}\,\right] 
&(\text{if }\ c\in(0,\frac{1}{2}))
\\[1em]
\left[\,-\sqrt{v\beta_c},\,\sqrt{v\beta_c}\,\right]
&(\text{if }\ c\in[\tfrac{1}{2},1)).
\end{cases}
\end{equation}
If $c=0$, then $\mu=\delta_0$.
If $c=1$, then $\mu$ is the semicircle law  on $[-2\sqrt{v}, 2\sqrt{v}]$.
\end{Theorem}

\begin{Remark}
{\rm
The formula \eqref{eq:mu} is valid for the extreme cases $c=0$ or $c=1$.
If $c=0$ then {there is no density} and $\mu=\delta_0$.
If $c=1$, then it can be checked that
$\sqrt[3]{R_+(x;\,1)}-\sqrt[3]{R_-(x;\,1)}=\sqrt{3}x\,\sqrt{4-x^2}$
so that, {for $v=1$}
we get  the semicircle law
$\mu(dt)= (1/2\pi) \sqrt{4-t^2}\chi_{[-2, 2]}(t)dt$
of ~\citet{Wigner55}.
}
\end{Remark}

\begin{proof}[Sketch of the proof]
We first derive the Stieltjes transform of the limiting eigenvalue distribution
by applying Theorem~\ref{theo:variance profile method}  to $Y_n=U_n/\sqrt{n}$.
Let $U_n=(U_{ij})_{1\le i,j\le n}$, so that $Y_{ij}=(1/\sqrt{n})U_{ij}$.
The variance profile is given as
\begin{equation}
    \label{def:profile}
\sigma(x,y)=    
\left\{
\begin{array}{l}
  v \quad \text{if}\ (x,y)\in \mathcal{C},\\
  0 \quad \text{otherwise,}
\end{array}
\right.
\quad
\mathcal{C}:= \set{(x,y)\in[0,1]^2}{\min(x,y)\le c}.
\end{equation}
We check easily that the conditions \eqref{eq:moment} are satisfied, since,   by   \eqref{def:gaussianU} and writing
$M:=\max\{|v-v'|, v',v  \}$, we get
\[
\delta_0(n)\le \frac{3M}{n}
\quad \textrm{and } \quad
\max_{i,j\le n} 
\frac{\mathbb{E}(Y_{ij}^4)}{n(\mathbb{E}Y_{ij}^2)^2}
\le \frac{\max\{{{M_4},{M_4}'}\}}{n\min\{v,v'\}}.
\]

Let us fix $z\in\C^+=\set{z\in\C}{\mathrm{Im}\,z>0}$.
The functional equation \eqref{eq:FE of vp}
from  Theorem~\ref{theo:variance profile method} becomes
\[
\eta_z(x)=
\begin{cases}
-\left(z+v\int_0^1\eta_z(y)\,dy\right)^{-1}&(x\le c),\\[0.8em]
-\left(z+v\int_0^c\eta_z(y)\,dy\right)^{-1}&(x>c).
\end{cases}
\]
Observe that the
right-hand sides are independent of $x$.
Integrating both sides of these equations,
we obtain the following simultaneous equations
\begin{equation}
\label{eq:FE1}
B=\frac{-c}{z+v A},\quad
A-B=\frac{c-1}{z+v B},
\end{equation}
where $A=\int_0^1\eta_z(x)\,dx$ and $B=\int_0^c\eta_z(x)\,dx$.
Note that $A$ is the desired Stieltjes transform $S(z)$.

If $c=0$, then we have $A=-1/z$ so that the 
limiting measure  is $\mu=\delta_0$.
If $c=1$ then 
the equation \eqref{eq:FE of vp}
reduces to the equation $A=-(z+vA)^{-1}$, which corresponds to the Stieltjes transform of the semi-circular law (cf.\ \citet[p.178]{Tao}).
Thus we assume $0<c<1$ in what follows.

Then,
the cubic equation for $A$, resulting from~\eqref{eq:FE1} writes
\begin{equation}
\label{cub3}
zA^3+(2z^2+1-2c)A^2+(z^2+2-2c)zA+z^2-c^2=0
\end{equation}
and it is an algebraic equation  with polynomial coefficients.
{The last
} equation \eqref{cub3} is reduced to
\begin{equation}
\label{eq:cubic1}
Y^3+p\left(z_v\right)Y+q\left(z_v\right)=0,
\end{equation}
where we set $z_v:={z}/{\sqrt{v}}$,
\begin{equation}\label{eq:YA}
Y=Y(z):=\frac{vA}{z}+\frac{2}{3}-\frac{(2c-1)v}{3z^2},
\end{equation}
and the coefficients
$p,q$ are given by the following analytical rational functions on $\C^*{:=\C\setminus\{0\}}$
\[
\begin{array}{l}
\ds p(z):=-\frac{z^4-2(c+1)z^2+(2c-1)^2}{3z^4},\\[1em]
\ds q(z):=-\frac{2}{27}\cdot\frac{z^6-3(c+1)z^4+\frac{3}{2}(5c^2-2c+2)z^2+(2c-1)^3}{z^6}.
\end{array}
\]
The  discriminant  of {the cubic equation} \eqref{eq:cubic1} is expressed by $p(z)$ and $q(z)$, { using $\alpha_c,\beta_c$ in \eqref{eq:alphabeta}}, as
(cf.\ \citet{Ronald})
\[
\Disc(z)=-\bigl(4p(z)^3+27q(z)^2\bigr)=\frac{4c^2(1-c)}{z^{10}}
(z^2-\alpha_c)(z^2-\beta_c).
\]
Let $\mathcal{E}=\set{z\in\C}{z=0\text{ or }\Disc(z_v)=0}$ be the set of exceptional points of \eqref{eq:cubic1}. For 
$z\not\in \mathcal{E}$, the equation \eqref{eq:cubic1} has three different solutions
(cf.\ \citet{Ronald}).
Cardano's method and formula
\eqref{eq:YA} imply  that, for $z\in\C^+$
\begin{equation}
\label{eq:StieltjesFromulaONE}
S(z)=\frac{z(u_+(z)+u_-(z))}{3v}
-\frac{2z}{3v}+\frac{2c-1}{3z}
\end{equation}
with $u_\pm(z):=\left(F_c(z_v)\pm\iunit \, D_c(z_v)\right)^{\frac13}$, $F_c(z):=-\frac{27}{2}q(z)$ and 
\[
D_c(z):=27\cdot \sqrt{\frac{\Disc(z)}{4\cdot 27}}=\frac{3c\sqrt{3-3c}}{z^5}\,
\sqrt{(z^2-\alpha_c)(z^2-\beta_c)},
\]
where  convenient branches of the cube and the square roots 
are chosen, respectively,  
for 
$u_\pm(z)$ and $D_c(z)$
to be such that $S(z)$
is a Stieltjes transform of a probability measure.
In particular, $S(z)$ is holomorphic on $\C^+$  and 
\begin{equation}
\label{eq:condition of StieltjesONE}
   u_+(z)\cdot u_-(z)=-3p(z),\quad
   \textrm{and} \quad\mathrm{Im}\,S(z)>0
\quad (z\in\C^+).
\end{equation}
Note that the branches of the roots    may be different on different subregions of $\C^+$
and that $U:=(u_++u_-)/3$
is a solution of \eqref{eq:cubic1}.
In order to derive the limiting eigenvalue distribution $\mu$ from $S(z)$,
we will need
the following  properties of  $S(z)$.
Set $\R^*:=\R\setminus\{0\}$.

\begin{Proposition} \label{prop:AnalCont}
 The limit $\ds S(x)=\lim_{y\to +0} S(x+yi)$  exists
for each $x\in\R^*$. The function 
 $S$ is continuous on $ \R^*$ and $S(x)$ is a solution
 of \eqref{cub3} on $ \R^*$.
\end{Proposition}
\begin{proof}[Sketch of the proof  of the proposition]
It is sufficient to prove it for a solution 
$U(z)$ of the reduced equation
\eqref{eq:cubic1}  on $\C^+$,
such that $U(z)$ is holomorphic on $\C^+$.
We apply Theorem X.3.7 of \citet{Palka}
to a convenient connected and simply connected domain $D$ avoiding the set $\mathcal{E}$. 
By the discussion of \citet[p.304]{Ahlfors},
$U$ has at most an ordinary algebraic singularity at a non-zero exceptional point,
so $U(z)$ is continuous on
$\R^*$.
\qed
\end{proof}

Without loss of generality, we suppose $v=1$.
We first assume that $x=0$.
The detailed local analysis of \eqref{eq:StieltjesFromulaONE}  
and \eqref{eq:condition of StieltjesONE}
that we omit here,  shows that
\begin{enumerate}
    \item[(Z1)] if $0<c<\frac12$, then $\ds\lim_{y\to+0}y\textrm{Im}\,S(yi)=  1-2c$, 
    so $\mu$ has an atom at $0$ with the mass $1-2c<1$,
    \item[(Z2)] if $c=\frac12$, then $\ds\lim_{y\to+0} \textrm{Im}\,S(yi)=+\infty,\ 
\lim_{y\to+0}y\textrm{Im}\,S(yi)=0$ so  $\mu$ does not have an atom at $0$,
    \item[(Z3)] if $\frac12<c<1$, then $\ds \lim_{y\to+0} \textrm{Im}\, S(yi)=c (2c-1)^{-1/2}=  \pi f_c(0)$, so  $\mu$ does not have an atom at $0$.
\end{enumerate}
Next we consider the case $x\ne 0$.
Combining the fact that
$S(z)$ is an odd function as a function on $\C\setminus\R$ 
by \eqref{eq:StieltjesFromulaONE}
and the property $S(\overline{z})=\overline{S(z)}$ of the Stieltjes transform,
we obtain $\mathrm{Im}\,S(-x+iy)=\mathrm{Im}\,S(x+iy)$ so that
$\mathrm{Im}\,S(-x)=\mathrm{Im}\,S(x)$.
Thus we can assume that $x>0$.

Suppose $\Disc(x)\ge 0$. Since the coefficients $p,q$ of \eqref{eq:cubic1}
are real on $\R^*$,  
the equation~\eqref{eq:cubic1} has only real solutions (cf.\ \citet{Ronald}).
Therefore, $S(x)$ is real so that
the density of $\mu$ vanishes at such points.

Next we assume that $\Disc(x)<0$.
By Proposition \ref{prop:AnalCont}, 
$S(x)$ is a solution of the cubic equation \eqref{cub3} and
$U(x)=(u_+(x)+u_-(x))/3$ is a solution of the reduced equation \eqref{eq:cubic1}.
In particular, 
the formulas \eqref{eq:StieltjesFromulaONE} and \eqref{eq:condition of StieltjesONE}
hold for $S(x)$, 
with convenient choices of branches of cubic roots and square roots.
Consequently, we have
\begin{center}
$\bigl\{F_c(x)+i D_c(x),F_c(x)-i D_c(x)\bigr\}=\bigl\{R'_+(x),\,R'_{-}(x)\bigr\}$
\end{center}
as a set,
where $R'_{\pm}(x):=R_{\pm}(x;\,c)/x^6\in\R$.
Let $\omega=e^{2\iunit \pi/3}$ denote the cube root of $1$ with positive imaginary part.
Then, 
\eqref{eq:StieltjesFromulaONE} yields that
the sum $u_{+}(x)+u_-(x)$ has the following form
\[
u_+(x)+u_{-}(x)=\omega^{k_+}\sqrt[3]{R'_+(x)}+\omega^{k_-}\sqrt[3]{R'_-(x)}
\quad\text{with}\quad
k_+,k_-\in\{0,1,2\}.
\]
By the first condition in~\eqref{eq:condition of StieltjesONE},
as $p(x)\in \R$,
we need to have $k_++k_-\equiv 0$ mod $3$,
that is, $(k_+,k_-)=(0,0)$, $(1,2)$ and $(2,1)$.
Using the fact that $R'_{+}(x)  > R'_{-}(x)$
when $x>0$ and $\Disc(x)<0$,
we see that the imaginary part
of $u_+(x)+u_-(x)$
and of $\lim_{y\to 0+} S(x+iy)$ is, respectively,
nul, positive and  negative
in these three cases. 
Since $\mathrm{Im}\, S(z)>0$, the last case is impossible.
Set $h(x):={\rm Im} \bigl(\omega\sqrt[3]{R'_+(x)}+ \omega^2\sqrt[3]{R'_-(x)}\bigr)$.
Notice that $h$ is a strictly positive continuous function on the set
$\set{x\in\R}{\Disc(x)<0}$ 
and
that $\frac{1}{\pi}h(t)=f_c(t)$, the density part of $\mu$ in the formula \eqref{eq:mu}.
Since the function $\mathrm{Im}\,S$ is continuous on $\R^*$ by Proposition~\ref{prop:AnalCont}, 
we have
$\mathrm{Im}\,S\equiv h$ or $\mathrm{Im}\,S\equiv 0$ 
on the set $\set{x\in\R^*}{\Disc(x)<0}$.

We now show that the latter case is impossible.
Note that 
$\mu$ has no atoms different from zero because $S(z)$ is continuous on $\overline{\C^+}\setminus\{0\}$.
By Theorem 2.4.3 of \citet{AGZ} and by the dominated convergence,
we have for closed intervals $[a,b]\subset\R^*$ 
\begin{equation}\label{eq:off1}
\mu([a,b])=
\frac1\pi  \lim_{y\to 0+}\int_a^b S(x+iy)\,dx
=\frac1\pi \int_a^b \lim_{y\to 0+} S(x+iy)\,dx=
0,
\end{equation}
so that $\mu(0,\infty)=0$ and, symmetrically, $\mu(-\infty,0)=0$.
Since $\mu$ is a probability measure, we get $\mu=\delta_0$.
This contradicts properties (Z1-3) proven in the case $x=0$.
Thus, we  have $\mathrm{Im}\,S\equiv h$ 
on the set $\set{x\in\R^*}{\Disc(x)\le 0}$
and, for $x\in\R^*$,
$\ds \lim_{y\to 0+} \frac1\pi \mathrm{Im}\,S(x+iy) =\frac{1}{\pi}h(x)=f_c(x)$.
Note that $f_c$ has a compact support $\{\Disc(x)\le 0\}$.
For $c\not=\frac12$, 
the function $f_c$ is continuous on $\R$. 
For $c=\frac12$, 
a detailed analysis shows that $\lim_{x\to 0}f_c(0)=\infty$,
with $f_c(x)\sim |x|^{-1/2}$ at $x=0$ and  $f_c$ is continuous on $\R^*$.
By property (Z3),
if $c>\frac12$
then $\lim_{y\to 0+} \mathrm{Im}\,S(iy) =\pi f_c(0)$.
When $c\not=1/2$, 
Proposition~\ref{th:0 ONE}.1 implies that $\mu=f_c(t)\,dt + [1-2c]_+\delta_0$.
Actually, 
if $s(z)$ is the Stieltjes transform of $\mu-f_c(t)\,dt - [1-2c]_+\delta_0$, then, using Proposition \ref{th:0 ONE}.2, we get
$\lim_{y\to 0+} {\rm Im}\, s(x+iy)=0$ for all $x\in \R$. 
When $c=1/2$,
by Proposition~\ref{th:0 ONE}.2, we get 
$\lim_{y\to 0+} {\rm Im}\, s(x+iy)=0$ for all $x\in \R^*$, uniformly on
compact intervals $[a,b]\subset \R^*$.
Like in \eqref{eq:off1}, we conclude
by Theorem 2.4.3 in \citet{AGZ} that $\mu=f_c(t)\,dt$.
The  support formula \eqref{eq:supp} follows by $\mathrm{supp}\,f_c=\{\Disc(x)\le 0\}$.
\qed
\end{proof}

In the Figures \ref{fig:1}--\ref{fig:5} we present graphical comparison between simulations for $n=4000$ and the limiting densities, when $c=1/5,2/5, 1/2, 3/5, 4/5$. 
%\sout{The graphs of $\alpha_c$ and $\beta_c$ as functions on $c$ are given in the Figure \ref{fig:6}.}
%\bl{It is mentioned in p. 9}

%\captionsetup{format=hang, font=bf, margin=70pt, name=fig}

\begin{Remark}
{\rm
We can also consider the class of generalized daisy graphs $D(a,b,k)$, 
with $b$ complete satellites  of $k$ vertices, so that there are  $N=a+kb$ vertices.  
If all three sequences $a_n,b_n,k_n$ are non-decreasing, the
graphs $D(a_n,b_n,k_n)$ form a growing sequence of graphical models.
Let us assume that  $k_n=k$ is fixed for $n$ large enough.
}
\end{Remark}
 \begin{Corollary}
 Consider a sequence of growing 
 graphs $D_n:=D(a_n,b_n,k)$ with $N_n=a_n+b_n k$ vertices. 
 Let $U_n$
 be a Wigner $N_n\times N_n$ matrix 
on $\mathcal{U}_{D_n}$ defined as in \eqref{def:gaussianU}. 
Assume that $\lim_{n\to+\infty}a_n/N_n=c\in[0, 1]$.
Then, the limiting eigenvalue distribution $\mu$ of the rescaled  matrices $U_n/\sqrt{N_n}$ exists
and is given by formula 
\eqref{eq:mu}.
\end{Corollary}
 \begin{proof}
The proof of Theorem \ref{theo:LEDforGW} is valid for the matrices  $U_n$ of size
 $N_n\times N_n$.
 Actually, the variance profiles $\sigma$ are the same and are  given by
 \eqref{def:profile} for all cases
$D(a_n,b_n,k)$. 
There are at most $(k^2+2)N_n$ non-zero perturbation terms $\delta_{ij}(N_n)$ and they are all bounded by $M=\max\{|v-v'|, v',v  \}$  so that
 $\delta_0(N_n)=O(1/N_n)\to 0$.
 \qed
 \end{proof}
\begin{Remark}
{\rm
The Wigner case may be considered in a framework of operator-valued free probability theory by methods of the rectangular free probability 
(cf.\ \citet[Chapter 9]{MingoSpeicher}, \citet{BenaychGeorges}).
}
\end{Remark}

\section{Wishart Ensembles of  Vinberg Matrices}
\label{sect:Wishart}

In this section,
we shall consider the quadratic Wishart (covariance) matrices introduced in \S\ref{ssect:WQ}.
We first prepare some special functions which we need later.
They generalize the Lambert $W$ function appearing (see \citet{Cheliotis})
in the case $\dVin=\mathrm{Sym}(n,\R)^+$ and $\ul{m}=(1,\dots,1)$.

\subsection{Lambert--Tsallis $W$ function and Lambert--Tsallis function $W_{\kappa,\gamma}$}
\label{sect:functionW}

For a non zero real number $\kappa$, we set
\[
\exp_\kappa(z):=\left(1+\frac{z}{\kappa}\right)^\kappa\quad(1+\frac{z}{\kappa}\in\C\setminus\R_{\le 0}),\quad
\tlog[\kappa](z):=\frac{z^\kappa-1}{\kappa}\quad(z\in\C\setminus\R_{\le 0}),
\]
where we take the main branch of the power function when $\kappa$ is not integer.
If $\kappa=\frac{1}{1-q}$, then it is exactly the so-called Tsallis \textit{$q$-exponential function} and \textit{$q$-logarithm}, respectively
(cf.\ \citet{AmariOhara2011,ZNS2018}).
We have the following relationship between these two functions:
\begin{equation}
    \label{rel:q exp log}
    \tlog[1/\kappa]\circ\exp_\kappa(z)=z\quad(-\pi<\kappa\mathrm{Arg}\left(1+\frac{z}{\kappa}\right)<\pi).
\end{equation}
By virtue of
$\ds\lim_{\kappa\to\infty}{\exp_\kappa(z)}=e^z$,
we regard $\exp_\infty(z)=e^z$ and $\tlog[0](z)=\log(z)$.

For two real numbers $\kappa,\gamma$ such that $\gamma\le\frac{1}{\kappa}\le 1$ and $\gamma<1$,
we introduce a holomorphic function 
$f_{\kappa,\gamma}(z)$,
which we call \textit{generalized Tsallis function}, by
\[
f_{\kappa,\gamma}(z):=\frac{z}{1+\gamma z}\exp_\kappa(z)\quad(1+\frac{z}{\kappa}\in\C\setminus\R_{\le 0}).
\]
{We note that $\kappa\in(-\infty,0)\cup [1,+\infty)$.}
Analogously to Tsallis $q$-exponential, 
we also consider $f_{\infty,\gamma}(z)=\frac{ze^z}{1+\gamma z}$ $(z\in\C)$.
In particular, $f_{\infty,0}(z)=ze^z$.
\out{We will write $\kappa\in[1,+\infty]$ if $\kappa\in[1,+\infty)$ or $\kappa=\infty$.}

In our work it is crucial to consider an inverse function to $f_{\kappa,\gamma} $.  
A multivariate inverse function of $f_{\infty,0}(z)=ze^z$ is called 
the Lambert $W$ function and studied
in \citet{Corless}. 
Hence, we call an inverse function to $f_{\kappa,\gamma}$ 
the {\it Lambert--Tsallis $W$ function.}

The function $f_{\kappa,\gamma}(z)$ has the inverse function {$w_{\kappa,\gamma}$} 
in a  neighborhood of $z=0$,
because we have $f'_{\kappa,\gamma}(0)=1\ne 0$ by
\[
f'_{\kappa,\gamma}(z)=\frac{\gamma z^2+\bigl(1+1/\kappa\bigr)z+1}{(1+\gamma z)^2}\left(1+\frac{z}{\kappa}\right)^{\kappa-1}.
\]

Let us present some properties of $f_{\kappa,\gamma}$.
When $\gamma\kappa\not=1$, the function $f_{\kappa,\gamma}$  
has a pole at $x=-\frac{1}{\gamma}$.
By the condition on $\kappa$ and $\gamma$,
the function $\gamma z^2+(1+1/\kappa)z+1$ has two real roots,
say $\alpha_1\le\alpha_2$, when $\gamma\not=0$.
If $\gamma=0$, there is only one real root, 
that we denote $\alpha_2=-\frac{\kappa}{\kappa+1}$.
$f_{\kappa,\gamma}'(z)=0$ implies $z=\alpha_i$  $(i=1,2)$, or $z=-\kappa$ if $\kappa>1$.
For the case $\kappa<0$,
it is convenient to change the variable by a homographic action $z'=\frac{z}{1+\frac{z}{\kappa}}$. Then
\[%\begin{equation}\label{kprime}
f_{\kappa,\gamma}(z)=f_{\kappa',\gamma'}(z')\quad\text{where}\quad
\kappa'=-\kappa>0,\quad
\gamma'=\gamma-\frac{1}{\kappa}.
\]%\end{equation}
Since a homographic action by an element in $SL(2,\R)$ leaves $\C^+$ invariant,
the analysis of the case $\kappa<0$ reduces
to the case $\kappa'>0$ and $\gamma'\le 0$.
 Then, the set $\mathcal{S}:=\R\setminus f_{\kappa,\gamma}(\R)$
has the following  possibilities.

\begin{enumerate}
    \item[(S1)] $\mathcal{S}=(f_{\kappa,\gamma}(\alpha_2),f_{\kappa,\gamma}(\alpha_1))$, where $f_{\kappa,\gamma}(\alpha_2)<f_{\kappa,\gamma}(\alpha_1)<0$.
    It occurs when 
    {$\kappa\in[1,+\infty]$} and $\gamma<0$, {and when $\kappa<0$ and $\gamma'=\gamma-\frac{1}{\kappa}<0$}.
    \item[(S2)] $\mathcal{S}=(-\infty,f_{\kappa,\gamma}(\alpha_2))$, where $f_{\kappa,\gamma}(\alpha_2)<0$.
    It occurs when 
      $\kappa>1$ and $\gamma\ge 0$ and when $(\kappa,\gamma)=(1,0)$.
    \item[(S3)] $\mathcal{S}=(-\infty,f_{\kappa,\gamma}(\alpha_1))$, where $f_{\kappa,\gamma}(\alpha_1)<0$.
    It occurs when 
    $\kappa<0$ and $\gamma'=\gamma-\frac{1}{\kappa}=0$.
    \item[(S{4})] $\mathcal{S}=(f_{\kappa,\gamma}(\alpha_1),f_{\kappa,\gamma}(\alpha_2))$, where $f_{\kappa,\gamma}(\alpha_1)<f_{\kappa,\gamma}(\alpha_2)<0$.
    It occurs when  $\kappa=1$ and $\gamma>0$.
\end{enumerate}
The cases (S1,S2,{S3}) are studied in detail in the Supplementary Material.
The case (S{4}) appears in the well known Wishart Ensemble case.

\begin{Theorem}\label{theo:LTfunction}
Let $\mathcal{S}$ be an interval or half-line given by {\rm(S1)-(S4)} above,
and $\overline{\mathcal{S}}\subset(-\infty,0)$ its closure.
Then, there exists a complex domain $\Omega\subset\C$,
symmetric with respect to the real axis and containing 0,
such that $f_{\kappa,\gamma}$ maps $\Omega$ bijectively to $\C\setminus\overline{\mathcal{S}}$.
Consequently,
the function $w_{\kappa,\gamma}$ can be continued in a unique way to 
a holomorphic function $W_{\kappa,\gamma}$ defined on 
$\C \setminus \overline{\mathcal{S}}$.
The codomain of $W_{\kappa,\gamma}$ is $\Omega$, {that is, $W_{\kappa,\gamma}(\C\setminus\overline{\mathcal{S}})=\Omega$}.\\
\end{Theorem}
\begin{proof}
The proof is based on the properties of $f_{\kappa,\gamma}$
showed in Proposition \ref{fkappagamma}.
\end{proof}

Recall that  the main branch of the Lambert {$W$} function is holomorphic 
on $\C\setminus (-\infty, -\frac1e]$
(see \citet{Corless}).

\begin{Definition}
{\rm
The unique holomorphic extension $W_{\kappa,\gamma}$
of $w_{\kappa,\gamma}$ to $\C \setminus \overline{\mathcal{S}}$
is called
{\it the main branch} of Lambert-Tsallis $W$ function.
In this paper, 
we only study and use $W_{\kappa,\gamma}$ among other branches
so that we call $W_{\kappa,\gamma}$ the \textit{Lambert--Tsallis function} for short.
Note that in our terminology the Lambert-Tsallis  $W$ function  is multivalued
and the
Lambert-Tsallis function $W_{\kappa,\gamma}$ is single-valued.
}
\end{Definition}

We summarize the basic properties of the Lambert-Tsallis function that we need later.

\begin{Proposition}\label{fkappagamma}
{\rm(i)}
Let $D=\Omega\cap\C^+$.
The function 
$f_{\kappa,\gamma}$
is continuous and injective on the closure
$\overline{D}$.
Consequently, $W_{\kappa,\gamma}$  extends continuously from $\C^+$ 
to $\C^+\cup\R$,
and one has $f_{\kappa,\gamma} (\partial\Omega\cap \C^+)= \mathcal{S}$.\\
%%%%%%%%%%%%%%%%%%%%%%%%%%%%%%%
\out{
\pg{and $W_{\kappa,\gamma}(\mathcal{S})=\partial\Omega\cap \C^+$
}.\\
\bl{It should be written as $\lim_{y\to+0}W_{\kappa,\gamma}(x+yi)\in\Omega\cap \C^+$ because $f_{\kappa,\gamma}(\partial \Omega\cap\la{\C^-})=\mathcal{S}$ (i.e. $W_{\kappa,\gamma}$ is not injective on $\mathcal{S}$).\\}
}
\noindent
{\rm(ii)}
The Lambert-Tsallis function $W_{\kappa,\gamma}$ 
has the following properties.
\begin{enumerate}
\item[\rm(a)] Suppose that $\kappa\ge 1$ and $\gamma<0$, or $\kappa<0$ and $\gamma'\le 0$.
In these cases, the set 
$D$ is bounded.
If $\kappa\ge 1$ then 
$D\subset \set{z\in\C^+}{\mathrm{Arg}\left(1+\frac{z}{\kappa}\right)\in(0,\frac{\pi}{\kappa+1})}$ and $z\in D$ satisfies $\mathrm{Re}\,z>-\kappa$.
{If $\kappa=\infty$, then 
%\sout{$\mathrm{Im}\,W_{\kappa,\gamma}(z)\in(0,\pi)$ for $z\in\C^+$}
{$D\subset\set{z\in\C^+}{\mathrm{Im}\,z\in(0,\pi)}$}.}
If $\kappa<0$ then  $D\subset \set{z\in\C^+}{\mathrm{Arg}\Bigl(\bigl(1+\frac{z}{\kappa}\bigr)^{-1}\Bigr)\in(0,\frac{\pi}{{|\kappa|+1}})}$.
Moreover,
$\lim_{|z|\to+\infty}W_{\kappa,\gamma}(z)=-\frac{1}{\gamma}$
(recall that $-\frac{1}{\gamma}$ is a pole of $f_{\kappa,\gamma}$).
\item[\rm(b)] Suppose  {$\kappa\in[1,+\infty]$} and $\gamma=0$. 
The set $D=\Omega\cap\C^+$ is unbounded and $f_{\kappa,0}(\infty)=\infty$.
If $\kappa\in[1,+\infty)$ then
$D\subset\set{z\in\C^+}{\mathrm{Arg}\left(1+\frac{z}{\kappa}\in(0,\frac{\pi}{\kappa+1})\right)}$.
If $\kappa=\infty$, then 
$W_{\infty,0}(z)$ is the classical Lambert function,
and one has 
%\sout{$\mathrm{Im}\,W_{\infty,0}(z)\in(0,\pi)$ for $z\in\C^+$}
{$D\subset\set{z\in\C^+}{\mathrm{Im}\,z\in(0,\pi)}$}.
\item[\rm(c)] Suppose $\gamma>0$. 
In this case we have $\kappa\in[1,\frac{1}{\gamma}]$.
The set $D=\Omega\cap\C^+$ is unbounded and $f_{\kappa,\gamma}(\infty)=\infty$.
Moreover, 
$D=\set{z\in\C^+}{\mathrm{Arg}\left(1+\frac{z}{\kappa}\right)\in(0,\frac{\pi}{\kappa})}$.
\end{enumerate}
\end{Proposition}

\begin{proof}
The main tool is the Argument Principle
(cf.\ \citet[Theorem 18, p.152]{Ahlfors}). 
A detailed study of the inverse image $f_{\kappa,\gamma}^{-1}(\R)$ is performed.
We omit the technical details, provided in Supplementary Material.
\qed
\end{proof}

\begin{Remark} 
{\rm
It is worth underlying that  we consider the main branch of the complex power function in the Tsallis {$q$-}exponential 
$\exp_\kappa(z)$ appearing inside the generalized Tsallis function $f_{\kappa,\gamma}$. Consequently,
the main branch  $W_{\kappa,\gamma}$ is the unique one such that $W(0)=0$.
A complete study of all branches of the Lambert-Tsallis  $W$ function will be interesting to do.
The study of the Lambert-Tsallis function $W_{\kappa,\gamma}$
in the full range of 
parameters $\kappa,\gamma$ is also an interesting open problem.
We exclude the case $\kappa\gamma>1$ with $\kappa>0$ because we do not need it later.
We note that, when
$\kappa\gamma>1$ 
and $\kappa>1$  with a condition $(1+\kappa)^2-4\gamma\kappa^2>0$,
then $f_{\kappa,\gamma}$ maps a subregion of $\C^+$ onto $\C^+$.
}
\end{Remark}

Applying the Lagrange inversion theorem, 
we see that the Taylor series
of the function $W_{\kappa,\gamma}$ near $z=0$ is
\[%\begin{equation}
\label{eq:taylorW}
W_{\kappa,\gamma}(z)=z+(\gamma-1)z^2+\left(\gamma^2-3\gamma+\frac{3\kappa+1}{\kappa}\right)z^3+o(z^3).
\]%\end{equation}

\subsection{Quadratic Wishart matrices}
\label{sect:4.2}

We will now study  eigenvalues of Wishart  (covariance) matrices in $P_n\subset 
\mathcal{U}_n$, defined in Section~\ref{ssect:WQ}.
We apply the approach of \citet[Cor.3.5]{LNofB},
based on the variance profile method (Theorem~\ref{theo:variance profile method}).

In this subsection, we first consider the case of $a_n=n-1$ and $b_n=1$, that is, $\dVin$ is  the symmetric  cone $\mathrm{Sym}(n,\R)^+$ of positive definite symmetric matrices of size $n$.
Let $\xi_n$ be a rectangular matrix of size $n\times N$.
In order to study eigenvalue distributions of $X_n=\xi_n{}^{\,t\!}\xi_n$,
we equivalently consider Wigner matrices of the form
\begin{equation}
\label{eq:mat of vp}
Y_n:=\pmat{0&\xi_n\\ {}^{t}\xi_n&0}\in\mathrm{Sym}(n+N,\,\R).
\end{equation}
If $X_n$ has eigenvalues $\lambda_j\ge0$ $(j=1,\dots,n)$,
then those of $Y_n$ are 
exactly $\pm\sqrt{\lambda_j}$ $(j=1,\dots,n)$ and zeros with multiplicity $|N-n|$.
Let $T_n$ denote the Stieltjes transform of the empirical eigenvalue distribution of rescaled $X_n/n$ and $S_n$  the Stieltjes transform  of rescaled $Y_n/\sqrt{n+N}$.
Then, it is easy to see that these Stieltjes transforms satisfy 
\begin{equation}
\label{formula:StieltjesT}
T_n\left(\frac{z^2}{p_n}\right)=
\frac{1}{2z}\left(\frac{1-2p_n}{z}+S_n(z)\right),
\end{equation}
where $p_n:=\frac{n}{n+N}$ and $q_n=\frac{N}{n+N}$.
In fact, we have \out{for $n\le N$}
\[
\begin{array}{r@{\ }c@{\ }l}
S_n(z)
&=&
% \ds\frac{1}{n+N}\left(\frac{|N-n|}{0-z}+\sum_{j=1}^{\min(n,N)}\Bigl(\frac{\sqrt{\lambda_j}}{\sqrt{n+N}}-z\Bigr)^{-1}+\Bigl(-\frac{\sqrt{\lambda_j}}{\sqrt{n+N}}-z\Bigr)^{-1}\right)\\
\ds\frac{1}{n+N}\left(\frac{|N-n|}{0-z}+\sum_{j=1}^{\min(n,N)}\frac{1}{\sqrt{\lambda_j}/\sqrt{n+N}-z}+\frac{1}{-\sqrt{\lambda_j}/\sqrt{n+N}-z}\right)\\
&=&
\ds
\frac{p_n-q_n}{z}
+
2zT_n\left(\frac{z^2}{p_n}\right).
\end{array}
\]

In order to study eigenvalue distributions of covariance matrices from Section~\ref{ssect:WQ}, with parameters $\ul{k}$ as in \eqref{def:k}, we introduce a trapezoidal variance profile $\sigma$ as follows.
Let $p,\alpha$ be real numbers such that 
$0< p< 1$ and $0\le \alpha\le (1-p)/p$.
Then, $\sigma$ is defined by
\begin{equation}
\label{def:varianec profile}
\sigma(x,y)=\begin{cases}
v&(x<p\text{ and }y\ge p+\alpha x),\\
v&(x\ge p\text{ and }0\le y\le \min\{(x-p)/\alpha,p\}),\\
0&(\text{otherwise}).
\end{cases}
\end{equation}
{Graphically, $\sigma$ is of the form
\begin{equation}
\label{eq:fig sigma}
\begin{array}{ccc}
    \sigma=\ \raisebox{-5em}{
    \includegraphics[scale=0.25]{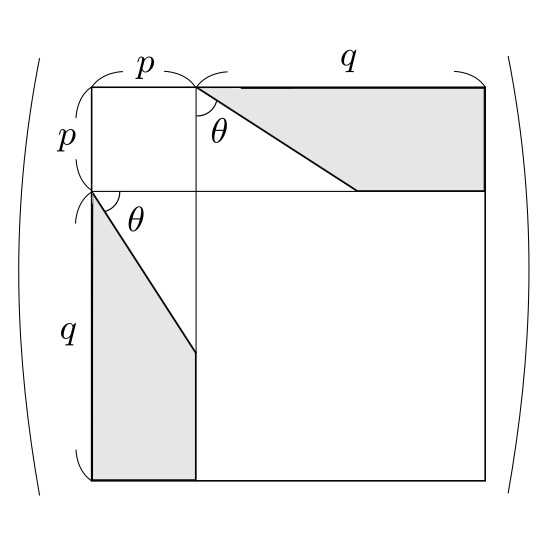}}&
    \text{with}&
    \begin{array}{l}
    p+q=1,\ p,q > 0\\
    0\le \tan \theta=\alpha\le \frac{q}{p}
    \end{array}
\end{array}
\end{equation}
}

If $\lim_n p_n=p$,
by Theorem~\ref{theo:variance profile method}, 
this variance profile determines the limiting distribution of
empirical eigenvalue distributions of the Wigner matrices $Y_n$ in~\eqref{eq:mat of vp}.
Recall that, to a variance profile $\sigma$,  
Theorem~\ref{theo:variance profile method} associates the Stieltjes transform $S_\sigma(z)$. 
It will be determined in Theorem \ref{theo:Stiltjes of vp}.
Analogously, to a variance profile $\sigma$ of  $\xi_n$, 
we associate
the ``covariance Stieltjes transform'' $T_\sigma(z)$ of
the corresponding covariance matrices ${Q_{\ul{k}}(\xi_n)=}\xi_n{}^{t}\xi_n$.
The covariance Stieltjes transform $T_\sigma(z)$ is
related to $S_\sigma(z)$ by the formula \eqref{formula:StieltjesT}. 
It will be determined in Proposition \ref{prop:Stieltjes of quad Sym}.

\begin{Theorem}
\label{theo:Stiltjes of vp}
Let $\sigma$ be a variance profile given in \eqref{def:varianec profile},
and set $\kappa:=1/(1-\alpha)$ and $\gamma:=(2p-1)/p=1-(q/p)$.
Then, the Stieltjes transform $S_\sigma(z)$ associated to $\sigma$ is given as
\[%\begin{equation}\label{eq:Sofsigma}
S_\sigma(z)=-\frac{2p}{zW_{\kappa,\gamma}\left(-\frac{vp}{z^2}\right)}+\frac{1-2p}{z}-\frac{2z}{v}
\quad(z\in\C^+),
\]%\end{equation}
where $W_{\kappa,\gamma}$ is the Lambert-Tsallis function defined in
Section \ref{sect:functionW}.
\end{Theorem}
\begin{proof} %(Sketch of the proof)
We use Theorem~\ref{theo:variance profile method}.
Let $z\in\C^+$ with $\mathrm{Im}\,z\gg 1$.
By \eqref{eq:FE of vp} we have
%By definition of $\sigma$ and $\eta_z$,
\[%\begin{equation}\label{eq:defetasup}
\eta_z(x)=\begin{cases}
\displaystyle
-\left(z+v\int_{p+\alpha x}^1\eta_z(y)\,dy\right)^{-1}
&
(0\le x\le p),\\
\displaystyle
-\left(z+v\int_0^{\alpha^{-1}(x-p)}\eta_z(y)\,dy\right)^{-1}
&
(p< x\le p+\alpha p),\\
\displaystyle
-\left(z+v\int_0^p\eta_z(y)\,dy\right)^{-1}
&
(p+\alpha p< x\le 1).
\end{cases}
\]%\end{equation}
For $z$ fixed, we set
\[
a(t):=\eta_z(t),\quad t\in[0,p],\qquad
b(t):=\eta_z(p+\alpha t), \quad t\in(0,p].
\]
By differentiating both sides in the above equations,
we obtain a  system
%differential equation
\begin{equation}
\label{def:difEqQW}
\left\{
\begin{array}{r@{\ }c@{\ }l}
a'(t)&=&-v\alpha a(t)^2b(t),\\
b'(t)&=&v a(t)b(t)^2,
\end{array}
\right.
\end{equation}
with 
initial data
$a(p)=-\left(z+v\int_{p+\alpha p}^1\eta_z(y)\,dy\right)^{-1}$, 
$b(0+)=-\frac{1}{z}$.
%By the property $\eta_z(x)\in\C^+$
By the unicity part of Theorem \ref{theo:variance profile method} holding  for $\eta_z(x)\in\C^+$,
it is enough to show that 
\eqref{def:difEqQW} is 
satisfied by 
\[
a(t)=
-zw(z)X(t)^{\alpha\kappa}
,\quad
b(t)=
-\frac{1}{z}\cdot X(t)^{-\kappa},
\]
where we set
$
w(z):=-\frac{1}{vp}W_{\kappa,\gamma}\left(-\frac{vp}{z^2}\right) \text{and}\ 
X(t):=1-\frac{vw(z)}{\kappa}\,t
$, 
and that $a(t),b(t) \in \C^+$ for Im\,$z$ big enough.
Here, we choose the main branches for complex power functions.
If $\alpha=1$ then
\[
a(t)=-zw(z)e^{-vw(z)t},\quad
b(t)=-\frac{1}{z}\cdot e^{vw(z)t}.
\]
The crucial part of the proof is to show that $a(t)$ satisfies the initial data condition.
We only give a proof for this in the case of $\alpha\ne 1$.
Set $w=w(z)$ and $X=X(p)$ for brevity.
Since $f_{\kappa,\gamma}(-vpw(z))=-\frac{vp}{z^2}$,
we have
\[
\begin{array}{l@{\quad\Longleftrightarrow\quad}l}
\ds
\frac{wX^\kappa}{1+v(1-2p)w}=\frac{1}{z^2}
&
\ds
wz^2X^\kappa=1+v(1-2p)w\\
&
\ds
wz^2X^\kappa=1-\frac{vwp}{\kappa}-(p+\alpha p-1)vw\\[1em]
&
\ds
X=z^2wX^\kappa+(p+\alpha p-1)vw\\
&
\ds
1=zwX^{\kappa-1}\left(z+(p+\alpha p-1)\frac{v}{z}\cdot X^{-\kappa}\right)\\
&
\ds
-zwX^{\kappa-1}=-\left(z+\frac{v(p+\alpha p-1)}{z}\cdot X^{-\kappa}\right)^{-1}.
\end{array}
\]
In the second and third equivalences, we use the formulas $\kappa=1/(1-\alpha)$ and $X=1-\frac{vwp}{\kappa}$.
Since $a(p)=-zwX^{\alpha\kappa}=-zwX^{\kappa-1}$ by $\alpha\kappa=\kappa-1$,
we see that
\[
a(p)=-\left(z+v\cdot\frac{p+\alpha p-1}{zX^\kappa}\right)^{-1}.
\]
Since $\eta_z(x)$ is independent of $x$ when $x\in[p+\alpha p,1]$, we have
\[
\int_{p+\alpha p}^1\eta_z(y)\,dy
=
(1-p-\alpha p)\eta_z(p+\alpha p)
=
(1-p-\alpha p)b(p)
=
\frac{p+\alpha p-1}{zX^\kappa}.
\]
Thus we conclude that $a(t)$ satisfies the initial condition. 
We omitted other details of the proof.
\qed

\end{proof}

\begin{Remark}
{\rm
We call the parameter $\kappa$ of Lambert-Tsallis functions the {\it angle parameter}
since it depends only on the angle of the trapeze in \eqref{eq:fig sigma}.
If $\kappa=1$, then we have $\alpha=0$ so that the trapeze reduces to a rectangle.
If $\alpha=q/p$, i.e.\ $\kappa=p/(p-q)=1/\gamma$,
then the trapeze reduces to a triangle.
On the other hand,
the parameter $\gamma=\frac{2p-1}{p}=1-C$ depends
directly on the  shape parameter $C=q/p$.
We call $\gamma$
the {\it  shape parameter} of the  Lambert-Tsallis function.
Note that the geometric condition $0\le \alpha\le \frac{p}{q}$ is equivalent to 
the condition $\frac{1}{\kappa}\ge \gamma$.
The formula $\gamma=1-\frac{q}{p}$ shows that $\gamma\in (-\infty,1)$.
We have
\[
\kappa\in [1,\tfrac{1}{\gamma}]\text{ if }0\le\gamma< 1,
\quad\text{and}\quad
\kappa\in [1,\infty] \cup (-\infty,\tfrac{1}{\gamma}] \text{ if }\gamma<0.
\]
}
\end{Remark}

Now we give the covariance Stieltjes transform $T_\sigma(z)$ for the profile $\sigma$.
\begin{Proposition}
\label{prop:Stieltjes of quad Sym}
Let $\sigma$ be a variance profile defined in~\eqref{def:varianec profile} with parameters $p$ and $\alpha$.
Set $\kappa:=\frac{1}{1-\alpha}$ and $\gamma:=\frac{2p-1}{p}=1-\frac{q}{p}$.
Then, the covariance Stieltjes transform $T_\sigma(z)$ corresponding to the profile $\sigma$ is described as 
\begin{equation}\label{Tsigma}
T_\sigma(z)
=T_{\kappa,\gamma}(z)
:=
-\frac{1}{v}-\frac{1}{zW_{\kappa,\gamma}\bigl(-\frac{v}{z}\bigr)}-\frac{\gamma}{z}
=\frac{\exp_\kappa\bigl(W_{\kappa,\gamma}(-v/z)\bigr)-1}{v}
\end{equation}
for $z\in\C^+$,
and its $R$-transform $R(z)$ is given as 
\[
R(z)=-\frac{1}{z}-\frac{v\gamma}{1-vz}-\frac{v}{(1-vz)\tlog[1/\kappa](1-vz)}
\quad
(1-vz\in\C\setminus\R_{\le 0}
).
\]
\end{Proposition} 
\begin{proof}%[Idea of the proof]
The first equality of the formula for $T_\sigma(z)$ is given by the formula \eqref{formula:StieltjesT}, and the second by the definition of the Lambert-Tsallis function.
To prove the formula of $R$-transforms, we use the fact that
$-\pi<\kappa\mathrm{Arg}\left(1+\frac{W(z)}{\kappa}\right)<\pi$ for any $z\in\C^+$
coming by Proposition \ref{fkappagamma} (ii) and we use relation~\eqref{rel:q exp log}.
\qed

\end{proof}

Recall that $\Omega$ denotes the codomain of $W_{\kappa,\gamma}$.
By Proposition \ref{fkappagamma}, 
for each $x\in\mathcal{S}$, 
there are exactly two solutions of $f_{\kappa,\gamma}(z)=x$ in $z\in\partial\Omega$, 
which are conjugate complex numbers,
denoted by $K_+(x)$, $K_-(x)$, such that $\mathrm{Im}\,K_+(x)>0$.
Recall that $\alpha_1\le\alpha_2$ are zeros of the function $\gamma z^2+(1+1/\kappa)z+1$.
Then, we have the following theorem.

\begin{Theorem}\label{th:WishartProfile}
Let $\sigma$ be a trapezoidal variance profile defined by \eqref{def:varianec profile}.
Let $\mu_\sigma$ be the probability measure corresponding to
the associated covariance Stieltjes transform  $T_\sigma$
given by \eqref{Tsigma}.
{
Then, 
the density function $d_\sigma$ of $\mu_\sigma$ is given as
\begin{equation}
    \label{eq:densityWishart}
d_\sigma(x)=\begin{cases}
\ds\frac{1}{2\pi xi}\left(\frac{1}{K_-(-\frac{v}{x})}-\frac{1}{K_+(-\frac{v}{x})}\right)&(\text{if }-\frac{v}{x}\in\mathcal{S}),
\\
0&(\text{if }-\frac{v}{x}\in\R\setminus\mathcal{S}).
\end{cases}
\end{equation}
Moreover, 
one has the following possibilities. 

%%%%%%%%%%%%%%%%%%%%%%%%%%%%%%%%
\out{for $\kappa>0$.
When $\kappa<0$, the same is valid by replacing $(\kappa,\gamma)$ by $(\kappa',\gamma')=(-\kappa,\gamma-\frac{1}{\kappa})$ 
 {according to \eqref{kprime}}.}
}
%%%%%%%%%%%%%%%%%%%%%%%%%%%%%%%%%
\begin{enumerate}
    \item In the case {$p<q$ and $\frac{q}{p}\ne\alpha$ (i.e.\ 
    {$\kappa\ge 1$} and $\gamma<0$, or $\kappa<0$ and $\gamma'<0$),}
    the measure $\mu_\sigma$ is absolutely continuous and its density $d_\sigma(x)$ is continuous on $\R$.
    In particular, $\mu_\sigma$ has no atoms.
    Its support is given as
    \out{(recall that $\alpha_1\le\alpha_2<0$ are zeros of the function $\gamma x^2+\bigl(1+1/\kappa\bigr)x+1$)
    }
   \begin{equation}\label{support1}
   \supp \mu_\sigma=\left[-\frac{v}{f_{\kappa,\gamma}(\alpha_2)},-\frac{v}{f_{\kappa,\gamma}(\alpha_1)}\right]
    =\left[\frac{v}{\alpha_2^2}\left(1+\frac{\alpha_2}{\kappa}\right)^{1-\kappa},\,\frac{v}{\alpha_1^2}\left(1+\frac{\alpha_1}{\kappa}\right)^{1-\kappa}\right].
    \end{equation}

    \item {In the case $p=q=\frac{1}{2}$ or $\frac{q}{p}=\alpha$ (i.e.\ {$\kappa\ge 1$} and $\gamma=0$, or $\kappa<0$ and $\gamma'=0$),}
    the measure $\mu_\sigma$ is absolutely continuous.
    Its density $d_\sigma$ is continuous on 
    $\R^*$ and $\lim_{x\to+0} d_\sigma(x)=+\infty$.
   In particular, $\mu_\sigma$ has no atoms.
   Let
   $\alpha_0:=\alpha_2$ if $\kappa\ge 1$ and $\alpha_0:=\alpha_1=-1$ if $\kappa<0$.
    The support of $\mu_\sigma$ is given as 
    \begin{equation}\label{S2}
    \supp\mu_\sigma=\left[0,-\frac{v}{f_{\kappa,\gamma}(\alpha_0)}\right]
    =
    \left[0,\frac{v}{\alpha_0^2}\left(1+\frac{\alpha_0}{\kappa}\right)^{1-\kappa}\right].
    \end{equation}
    When $\kappa=\infty$, the measure $\mu_\sigma$
    is the Dykema-Haagerup measure $\chi_v$
    with support $[0,ve]$.
    
    \item In the case {$p>q$ (i.e.\ {$\kappa\ge 1$} and $0<\gamma<1$),}
    we have 
    $\mu_\sigma= d_\sigma(x) dx + {(1-\frac{q}{p})}\out{(p-q)}\delta_0$.
    The measure $\mu_\sigma$ has an atom at $x=0$ with mass ${1-\frac{q}{p}}\out{p-q}$.
    Recall that $\kappa\in[1,1/\gamma]$.
    When $\kappa>1$, the support of $\mu_\sigma$ is given by \eqref{S2}.
    The function $d_\sigma$ is continuous on 
    $\R^*$ and  $\lim_{x\to+0} d_\sigma(x)=+\infty$.  
    For $\kappa=1$ and  $-\infty<\gamma<1$, 
    the measure $\mu_\sigma$ is the Marchenko-Pastur law $\mu_C$ with parameter $C=\frac{q}{p}= 1-\gamma\in(0,1)$ and 
    $\supp d_\sigma=\left[v(1-\sqrt{C})^2,\,v(1+\sqrt{C})^2\right]$.
\end{enumerate}
\end{Theorem}
\begin{proof}
We use 
%the formula of $T_\sigma(z)$
%from
Proposition~\ref{prop:Stieltjes of quad Sym}.
Let $z=x+yi$. 
{By Proposition \ref{fkappagamma} (i)}
and the fact that $W_{\kappa,\gamma}(z)=0$ only if $z=0$, 
we see that
$l(x):=\lim_{y\to+0}\mathrm{Im} \,T_\sigma(x+iy)$ exists when $x\ne 0$
and that $l(x)=0$ when $-v/x\not\in\mathcal{S}$.

Assume that $x\ne 0$ and $-v/x\in \mathcal{S}$.
Set $a(x)+ib(x):=\lim_{y\to0+}W_{\kappa,\gamma}(-v/z)$. 
Since the function $f_{\kappa, \gamma}$ is continuous and injective on the closure  $\overline{D}\subset \overline{\C^+}$, the function
$a+ib$ is continuous. 
By Proposition \ref{fkappagamma} (i),
we have $b(x)> 0$ and {$a(x)+ib(x)=K_+(-\frac{v}{x})$}.
Since $\overline{\mathcal{S}}\subset (-\infty,0)$ by Theorem~\ref{theo:LTfunction},
%\out{Recall that if $s\in\mathcal{S}$ then $s<0$; thus} 
we have $-v/x<0$, that is, $x>0$.
Thus, we obtain for $-v/x\in\mathcal{S}$ with $x\ne 0$
\begin{equation}\label{lx}
\begin{array}{r@{\ }c@{\ }l}
l(x)
&=&
\ds
\lim_{y\to0+}\mathrm{Im}\,T_\sigma(x+yi)
=
\mathrm{Im}\left(-\frac{1}{v}-\frac{1}{x(a(x)+ib(x))}-\frac{\gamma}{x}\right)\\[1em]
&=&
\ds
-\frac{1}{2xi}\left(\frac{1}{K_+(-\frac{v}{x})}-\frac{1}{K_-(-\frac{v}{x})}\right)
=
\frac{b(x)}{x(a(x)^2+b(x)^2)}>0,
\end{array}
\end{equation}
and thus $l(x)$ is a continuous function on $\R^*$.
Therefore,
$x\in\R^*$ is included in the support of $\mu_\sigma$ if and only if $-v/x\in\overline{\mathcal{S}}$.
By \eqref{eq:StContfONE}, 
we have $d_\sigma(x)=\frac{1}{\pi}l(x)$, so that we obtain \eqref{eq:densityWishart}.

Let us consider the case (S1).
In this case, 
since $\mathcal{S}=(f(\alpha_2),f(\alpha_1))$ and $f(\alpha_1)<0$, 
we have
\[
x\in\supp\mu
\iff
f(\alpha_2)\le -\frac{v}{x}\le f(\alpha_1)<0
\iff
-\frac{v}{f(\alpha_2)}\le x\le -\frac{v}{f(\alpha_1)}.
\]
Recall that $\alpha_i$, $i=1,2$ are the {real} solutions of the equation $\gamma z^2+(1+1/\kappa)z+1=0$.
For a solution $\alpha$ of this equation,
we have by $1+\alpha/\kappa=-\alpha(1+\gamma\alpha)$
\[
f_{\kappa,\gamma}(\alpha)=\frac{\alpha}{1+\gamma\alpha}\left(1+\frac{\alpha}{\kappa}\right)^{\kappa}
=
-\alpha^2\left(1+\frac{\alpha}{\kappa}\right)^{\kappa-1},
\]
so that we arrive at the assertion 1.\ of the theorem.
The argument  for other two cases  is similar, and hence we omit it.

Next we consider the case $x=0$.
We present the case $\kappa\in[1,+\infty)$ and $\gamma=0$.
For $z\in \C^+$,
let us set $re^{i\theta}=1+\frac{W_{\kappa,\gamma}(-v/z)}{\kappa}$
$(r>0,\ \theta\in(0,\pi))$.
{By Proposition \ref{fkappagamma} (ii-b)},
the set $D=\Omega\cap\C^+$ is unbounded and $f_{\kappa,\gamma}(\infty)=\infty$.
Consequently,
if $z\to0$ in $\C^+$, or equivalently $-v/z\to\infty$ in $\C^+$, 
then we have $W_{\kappa,0}(-v/z)\to\infty$ and $r\to+\infty$.
Again by {Proposition \ref{fkappagamma} (ii-b)},
we see that $\theta\in(0,\frac{\pi}{\kappa+1})$ so that
$\sin\kappa\theta>0$ when $z=-v/(iy)\in\C^+$,
and thus
\[
\begin{array}{r@{\ }c@{\ }l}
\mathrm{Im}\,T(z)
&=&
\ds
\mathrm{Im}\,\frac{\exp_\kappa\Bigl(W_{\kappa,\gamma}(-v/z)\Bigr)-1}{v}
=
\mathrm{Im}\,\frac{(re^{i\theta})^\kappa-1}{v}\\[1em]
&=&
\ds
\mathrm{Im}\,\frac{r^\kappa\cos\kappa\theta-1+i r^\kappa\sin\kappa\theta}{v}
=
\frac{r^\kappa\sin\kappa \theta}{v}\to+\infty\quad(y\to+0).
\end{array}
\]
On the other hand,
$\mu_\sigma$ does not have an atom at $x=0$
because we have by $W_{\kappa,0}(-v/z)\to\infty$ and by $\gamma=0$
\[
yT(iy)=-\frac{y}{v}-\frac{1}{iW_{\kappa,\gamma}(-v/(yi))}-\frac{\gamma}{i}\to \gamma i =0 \quad(y\to+0).
\]
The proofs for other cases are similar, and hence we omit them.

The absolute continuity of $\mu_\sigma$ follows from Proposition \ref{th:0 ONE}, 
by considering $\mu_0:=\mu_\sigma - d_\sigma(x)dx$, 
or, 
in the case with atom at $x=0$, 
of $\mu_0:=\mu_\sigma - d_\sigma(x)dx - \gamma\delta_0 $ 
and using the fact that the Stieltjes transform $S_0(z)$ of $\mu_0$ 
satisfies $\lim_{y\to 0+} \mathrm{Im}\,S_0(x+iy)=0$ for all $x\in\R$.
The argument is similar as in the proof of Theorem~\ref{theo:LEDforGW}.
\qed
\end{proof}

In the following corollary, we give a real implicit equation for the density $d_\sigma$ analogous
to the Dykema-Haagerup equation \eqref{dykh}.
To do so, 
we introduce the following notation
\[
e_\kappa(z):=\left|\exp_\kappa(z)\right|\ge 0,\quad
\theta_\kappa(z)=\kappa\mathrm{Arg}\left(1+\frac{z}{\kappa}\right)\quad(z\in\C^+).
\]
If $\kappa=\infty$, we set $e_\kappa(z):=e^{\mathrm{Re}\,z}$ and $\theta_\kappa(z):=\mathrm{Im}\,z$.
Then, we have $\exp_\kappa(z)=e_\kappa(z)\bigl(\cos\bigl(\theta_\kappa(z)\bigr)+i\sin\bigl(\theta_\kappa(z)\bigr)\bigr)$.
\begin{Corollary}
\label{cor:dykh type density}
{\rm(i)} Suppose $v=1$ for simplicity.
For two real numbers $\kappa,\gamma$ such that $\gamma\le \frac{1}{\kappa}\le1$ and $\gamma<1$, 
the density $d_\sigma$ of the limiting law $\mu_\sigma$ satisfies the equation
\begin{equation}
\label{eq:cor dykh}
\begin{array}{r}
\ds
d_\sigma\left(\frac{\sin\bigl(\theta_\kappa(z)\bigr)}{b}\Bigl(1+\gamma a-\gamma b\cot\bigl(\theta_\kappa(z)\bigr)\Bigr)\bigl(e_\kappa(z)\bigr)^{-1}\right)\qquad\\[1em]
\ds
=
\frac{1}{\pi}\cdot e_\kappa(z)\sin\bigl(\theta_\kappa(z))
\end{array}
\end{equation}
for $z=a+bi\in\partial \Omega\cap\C^+$.
In particular, 
when $(\kappa,\gamma)=(\infty,0)$, the density $d_\sigma$ satisfies the equation \eqref{dykh} with $b=x$ and $a=-x\cot x$ $(x\in[0,\pi))$.\\
{\rm(ii)} 
If $\kappa\in[1,\infty]$ and $\gamma<0$,
then the correspondence $a\mapsto b=b(a)$ is unique
for each $z=a+bi\in\partial \Omega\cap\C^+$.  
Then, $a\in[\alpha_1,\alpha_2]$.
The same is true for $\kappa=\infty$ 
and $\gamma=0$ with $a\in [-1,+\infty)$.
\end{Corollary}
\begin{proof}
{
(i) Let $z=a+bi\in\partial D\cap\C^+$.
Then, it satisfies $f_{\kappa,\gamma}(z)\in\mathcal{S}$.
Suppose $f_{\kappa,\gamma}(z)=-\frac{1}{x}$,
and set $X=a+\gamma a^2+\gamma b^2$ and $Y=|1+\gamma z|^2=(1+\gamma a)^2+(\gamma b)^2$.
Notice that $X^2+b^2=(a^2+b^2)Y$.
The equation $f_{\kappa,\gamma}(z)=-\frac{1}{x}$ means that
\begin{eqnarray}
\label{eq:im zero 1}
\frac{e_\kappa(z)}{Y}\Bigl(X\cos\bigl(\theta_\kappa(z)\bigr)-b\sin\bigl(\theta_\kappa(z)\bigr)\Bigr)&=&-\frac{1}{x},\\
\label{eq:im zero 2}
X\sin\bigl(\theta_\kappa(z)\bigr)+b\cos\bigl(\theta_\kappa(z)\bigr)&=&0.
\end{eqnarray}
The latter one \eqref{eq:im zero 2} yields that $\cos\bigl(\theta_\kappa(z)\bigr)=-\frac{\sin\bigl(\theta_\kappa(z)\bigr)}{b}X$ so that
\[
-\frac{1}{x}=-\frac{e_\kappa(z)}{Y}\cdot\frac{\sin\bigl(\theta_\kappa(z)\bigr)}{b}(X^2+b^2)
\iff
\frac{1}{x}\cdot\frac{b}{a^2+b^2}=e_\kappa(z)\sin\bigl(\theta_\kappa(z)\bigr).
\]
On the other hand,
\eqref{eq:im zero 2} can be written as $X=-b\cot\bigl(\theta_\kappa(z)\bigr)$, and using this expression together with \eqref{eq:im zero 1}, we obtain 
\[
-\frac{1}{x}=\frac{e_\kappa(z)}{Y}\Bigl(
{-b\cot\bigl(\theta_\kappa(z)\bigr)}
\cos\bigl(\theta_\kappa(z)\bigr)-b\sin \bigl(\theta_\kappa(z)\bigr)\Bigr)=-\frac{b}{\sin\bigl(\theta_\kappa(z)\bigr)}\cdot \frac{e_\kappa(z)}{Y}
\]
and hence
\[
x=\frac{\sin\bigl(\theta_\kappa(z)\bigr)}{b}\cdot Y \bigl(e_\kappa(z)\bigr)^{-1}.
\]
It is easy to check that we have $Y=1+\gamma a+\gamma X$.
By \eqref{lx}, the density can be described as $d_\sigma(x)=\frac{1}{\pi x}\cdot\frac{b}{a^2+b^2}$ so that we obtain the formula \eqref{eq:cor dykh}.\\
(ii)
%\la{TODO:replace $x,y$ by $a,b$ YES IT WILL BE BETTER}
We shall show the part (ii) for $\kappa\in(1,\infty)$ and $\gamma<0$.
The other cases can be done by a similar way.
Let $z=a+bi\in D=\Omega\cap\C^+$.
Set $\theta(a,b)=\mathrm{Arctan}\frac{b}{\kappa+a}$ for $a>-\kappa$ and $b>0$.
By Proposition~\ref{fkappagamma} (ii-a),
we see that $\mathrm{Re}\bigl(1+\frac{z}{\kappa}\bigr)=1+\frac{a}{\kappa}>0$
and hence $\theta_\kappa(a+ib)=\kappa\theta(a,b)$.
Note that $\frac{\partial}{ \partial b}\theta_\kappa(a+ib)=\kappa\cdot\frac{\kappa+a}{(\kappa+a)^2+b^2}$.
For given $a>-\kappa$, 
set $g(y;\,a):=y\cot(\theta_{\kappa}(a+iy))$. 
Let $y_0>0$ satisfy $\theta(a,y_0)=\frac{\pi}{\kappa+1}$.
Then, we can show that $g(y;\,a)$ is monotonic decreasing for $y\in(0,y_0)$.

Set $h(y)=h(y;\,a):=a+\gamma a^2+\gamma y^2+g(y)$  
for the fixed $a>-\kappa$.
Recall that  $h(y;\,a)=0$ if and only if $z=a+iy\in \partial D\cap \C^+$.
As $\gamma<0$, we see that
the function $h(y):=a+\gamma a^2+\gamma y^2+g(y)$ is decreasing on $y\in(0,y_0)$
for each fixed $a>-\kappa$.
Since $\cot(\theta_\kappa(a+iy_0))=-\frac{\kappa+a}{y_0}$,
we see that $h(y_0;\,a)<0$.
By the fact that $\lim_{y\to+0}g(y;\,a)=1+\frac{a}{\kappa}$, 
we have $\lim_{y\to+0}h(y;\,a)=\gamma(a-\alpha_1)(a-\alpha_2)$.
Since $h$ is monotonic decreasing on $y\in(0,y_0)$,
if $a\in (\alpha_1,\alpha_2)$ then $\lim_{y\to+0}h(y;\,a)>0$ so that 
there exists a unique solution $y=b$ of $h(y;\,a)=0$ in $y\in(0,y_0)$ for each $a\in(\alpha_1,\alpha_2)$
by the intermediate value theorem,
whereas
if $\lim_{y\to+0}h(y;\,a)\le 0$ then there is no solution of $h(y;\,a)=0$ in $y\in(0,y_0)$.
Thus the correspondence $a\mapsto b=b(a)$ is unique for each $z=a+bi\in\partial \Omega \cap\C^+$.
\qed
}

\end{proof}

\begin{figure}
    \centering
    \includegraphics[scale=0.5]{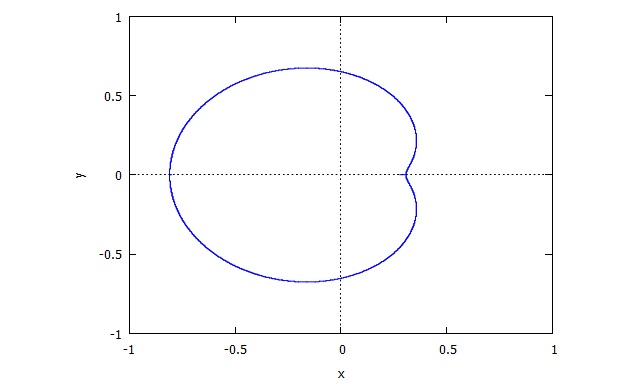}
     \caption{Domain $\Omega$ for $\kappa=-1/3, \gamma=-4$.}
     \label{fig:counterex}
\end{figure}

\begin{Remark}
{\rm
Corollary \ref{cor:dykh type density} (ii) enables us to write the density $d_\sigma$ with one real parameter in a way similar to \citet[Theorem 8.9]{DykemaHaagerup}, see formula \eqref{dykh}.
In particular, in the case (a), we obtain the formula
\[
d_\sigma\left(\frac{\sin b(a)}{b(a)}\Bigl(1+\gamma a-\gamma  b(a)\cot b(a)\Bigr)\,e^{-a}\right)
=
\frac{1}{\pi}\cdot e^a\sin b(a)
\quad(a\in[\alpha_1,\alpha_2]).
\]
}
\end{Remark}
%On conjecture: HN noticed that
A natural conjecture that we always have a $1$-$1$ correspondence $a\to b$ or $b\to a$ is not confirmed by numerical generation of the domain $\Omega$. For
 $\kappa=-1/3$ and $\gamma=-4$  the domain $\Omega$ is illustrated in the  Figure \ref{fig:counterex}.
We do not have unicity of $a\to b$ nor $b\to a$.

%\la{However, it is not true in general.
%It arises when $\kappa$ and $ \gamma$ are both negative, for example, $(\kappa,\gamma)=(-1/3,-4)$.}

\subsection{Applications to Wishart Ensembles of Vinberg matrices.}
\label{sect:4.3}

Now we apply Theorem \ref{th:WishartProfile}
to the covariance matrix $X_n=Q_{\ul{k}}(\xi_n)\in P_n$
in two situations.
The first (Corollary~\ref{prop:support}) is the case when $P_n$ is the symmetric cone $\mathrm{Sym}(n,\R)^+$ with $\ul{k}$ of the form \eqref{def:k type 1} below.
The second situation (Theorem~\ref{mainWishart}) is the  general case when $P_n\subset {\mathcal U}_n$ is a {dual} Vinberg cone with $\ul{k}$ of the form \eqref{def:k}.
This case contains the first one, that we present separately  because of the importance of the symmetric cone $\mathrm{Sym}(n,\R)^+$.

%\pg{WE COULD NOT PRESENT THE FIRST CASE IF WE STICK TO BREVITY.}

Let us assume that $\ul{k}=\ul{k}(n){=(k_1,\dots,k_n)}$ in~\eqref{def:k} is of the form
\begin{equation}
    \label{def:k type 1}
    \ul{k}=m_1(1,\dots,1,1)+m_2(n)(0,\dots,0,1),\quad
    \lim_n \frac{m_2(n)}{n}=m,
\end{equation}
where $m_1\in\mathbb{Z}_{\ge0}$ is a fixed non-negative integer 
and $m\in\R_{\ge0}$ is a non-negative real such that
$m_1+m>0$.
Set $N:=k_1+\cdots+k_n=m_1n+m_2(n)$.
We note that the case $m_1=0$ corresponds to the classical Wishart ensembles,
and  if $m_1\ge 1$ then we have $N\ge n$.

\out{
\bl{HN proposes $m_2=M_n$ such that $M_n/n$ converges as $n\to\infty$.}
\pg{I PROPOSE TO WRITE $m_2(n)$ INSTEAD OF $M_n$.}
\bl{HN agrees with it.}
}

\begin{Corollary}
\label{prop:support}
Let $\ul{k}$ be as in \eqref{def:k type 1}.
Suppose that $ \xi_n\in E_{\ul{k}}$ 
is an i.i.d. matrix with finite fourth moments and let $X_n=\xi_n\,{}^t\xi_n$.
Let $\mu_n$ be the empirical eigenvalue distribution of $X_n/n$. 
Then, there exists a limiting eigenvalue distribution $\mu=\lim_n \mu_n$.
The  Stieltjes transform $T(z)$ of $\mu$ is given by
formula \eqref{Tsigma}
\[
T(z)
=
T_{\kappa,\gamma}(z)
=
\frac{\exp_\kappa\Bigl(W_{\kappa,\gamma}(-v/z)\Bigr)-1}{v}
\ \text{with}\ 
\kappa=\frac{1}{1-m_1},\ 
\gamma=1-m-m_1.
\]
The measure $\mu$ is absolutely continuous and has no atoms.
If $m_1=0$ then the measure $\mu$ is the Marchenko-Pastur law with parameter $C=m$.
The case $(m_1,m)=(1,0)$ corresponds to the Dykema--Haagerup measure $\chi_v$.
If $m=0$ then the density $d$ is continuous on $\R^*$ and $\lim_{x\to+0}d(x)=+\infty$.
When $m_1\ge 2$ then 
the support of $\mu$ is $[0,vm_1^{m_1/(m_1-1)}]$.
Otherwise, for $m_1,m>0$,
the density $d(x)$ of $\mu$ is continuous on $\R$, and its support 
equals  $[A(\alpha_2),A(\alpha_1)]$ where
{$A(\alpha_i):= v\alpha_i^{-2}(1+(1-m_1)\alpha_i)^{m_1/(m_1-1)}$}
and $\alpha_1<\alpha_2$  are roots of the function
$(1-m_1-m)x^2+(2-m_1)x+1$.

%%%%%%%%%%%%%%%%%%%%%%%

\end{Corollary}
\begin{proof}
We use Theorem \ref{theo:variance profile method}.
It is enough to show that the matrix $Y_n$ in \eqref{eq:mat of vp} has the variance profile $\sigma$ in \eqref{def:varianec profile}
and that  the conditions \eqref{eq:moment} 
are satisfied.
Since we have for $n$ large enough
\[\bigl|\delta_0(n)\bigr|\le \frac{1}{n^2}\cdot 2v(m_1+m+1)n=\frac{2v(m_1+m+1)}{n}\to0\quad(n\to\infty)\]
and
if $\mathbb{E}|Y_{ij}|^2\ne0$ then 
\[\frac{\mathbb{E}(Y_{ij}^4)}{n(\mathbb{E}Y_{ij}^2)^2}=\frac{{M_4}}{vn}\to0\quad(n\to\infty),\]
we can easily check the conditions \eqref{eq:moment}.
%%%%%%%%%%%%%%%%%%%%%%%%%
Thus, we can apply Theorem
\ref{th:WishartProfile}.
Consider $m_1\ge 2$. Then 
$\kappa<0$.
When $m=0$, then we have $\gamma'=\gamma-\frac{1}{\kappa}=0$ so that we apply Theorem~\ref{th:WishartProfile}.2.
We have $\alpha=-1$, $1-\frac{1}{\kappa}=m_1$ and $1-\kappa=\frac{m_1}{m_1-1}$.
By \eqref{S2}, the support is given by $\supp \mu=\Bigl[0,\frac{v}{\alpha^2}\Bigl(1+\frac{\alpha}{\kappa}\Bigr)^{1-\kappa}\Bigr]
=
%\Bigl[0,vm_1^{\frac{m_1}{m_1-1}}\Bigr]
\Bigl[0,vm_1^{m_1/(m_1-1)}\Bigr]$.
When $m>0$, we have $\gamma'<0$ so that we apply Theorem~\ref{th:WishartProfile}.1.
The support of $\mu$
is given by the formula~\eqref{support1}, where
 $\alpha_1\le \alpha_2$ are roots of the function
$ \gamma x^2+
(1+1/\kappa)x+1$.
\qed
\end{proof}

\begin{Remark}
{\rm
If $m=0$,
our results contain those of
\citet[Section 4.5.1]{CR} and \citet[Theorem 4 and (12)]{Cheliotis}.
The result on the limiting densities of biorthogonal ensembles 
in \citet{Cheliotis} can be reproduced 
from Corollary~\ref{prop:support}.
In fact, our random matrices $Q_{\ul{k}}(\xi_n)$
%\out{$Q_{\ul{k}}(\xi_n{}^{\,t\!}\xi_n)$} 
essentially correspond to those considered in \citet{Cheliotis} 
through adjusting parameters
 $m_1=\theta-1$ and $m_2(n)=b-1$ (not depending on $n$), 
where $\theta$ and $b$ are parameters used in that paper.
}
\end{Remark}

\begin{figure}
    \centering
    \includegraphics[scale=0.3]{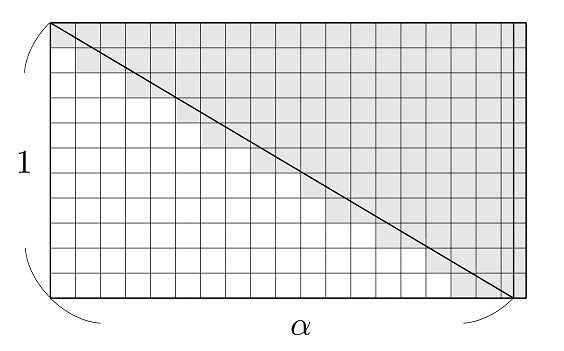}
    \caption{ Realization of non-integer $\alpha$}
    \label{fig:notN}
\end{figure}

\begin{Remark}
{\rm
Until now,
we assumed that $m_1\in\mathbb{Z}_{\ge0}$ and hence the parameter $\alpha$
of the variance profile $\sigma$
needs to be also an integer.
However,
we can take a sequence $\{\ul{k}(n)\}_{n=1}^{\infty}$ so that the corresponding $\alpha$ is an arbitrary given positive real number.
In fact, when $\alpha>0$ is given, 
we consider a right triangle with lengths $1$ and $\alpha$.
For an arbitrary $n$,
we cover the triangle by $1/n \times 1/n$ squares as in the figure.
To each $j=1,\dots,n$,
we associate an integer $k_j(n)$ such that $\frac{k_j(n)}{n}\le\frac{j}{n}\alpha< \frac{k_j(n)+1}{n}$, or equivalently $k_j(n)\le j\alpha<k_j(n)+1$,
and we set $k(n)=(k_1(n),\dots,k_n(n))$.
Note that this condition is independent of $n$ so that $k_j(m)=k_j(n)$ when $m\ge n\ge j$,
and hence $\{E_{\ul{k}(n)}\}_n$ is a sequence of vector spaces such that $E_{\ul{k}(n)}\subset E_{\ul{k}(n+1)}$.
In the Figure \ref{fig:notN}, we set $\alpha=1.8$, $n=11$ and $k(n)=(1,2,2,1,2,2,1,2,2,1,2)$.
}
\end{Remark}

Let us return to the quadratic Wishart case for general $\dVin$ with parameter $\ul{k}$ as in~\eqref{def:k} such that $m_1,m_2 \in\mathbb{Z}_{\ge0}$ are fixed.
Note that $m_2(n)$ in the previous discussion is now $m_2(n)=m_2b_n$.
Set $N_n:=m_1n+m_2b_n$.
We have
\[
E_{\ul{k}}=\set{\xi=\pmat{\eta\\ \zeta}\in\mathrm{Mat}(n\times N_n,\,\R)}{
\begin{array}{l}
\eta=(\eta_{ij})\in\mathrm{Mat}(a_n\times N_n,\,\R),\\ 
\zeta=(\zeta_{ij})\in\mathrm{Mat}(b_n\times N_n,\,\R)\\
\eta_{ij}=0\text{ if }j\le (m_1-1)i,\\
\zeta_{ij}=0\text{ if }M(i,j)\not\in\{1,2,\dots,m_1+m_2\}
\end{array}
},
\]
where $M(i,j):=j-m_1a_n-(m_1+m_2)(i-1)$.

\begin{Corollary}\label{mainWishart}
Let $\{\dVin\}_n$ be a sequence of generalized dual Vinberg cones such that $\lim_{n\to\infty}a_n/n=c\in(0,1]$.
Let $\ul{k}$ be a vector as in~\eqref{def:k} such that $m_1,m_2$ are fixed.
Set $\kappa:=1/(1-m_1)$ and $\gamma:=1-\bigl(m_1+m_2(1-c)\bigr)/c$.
Then, the Stieltjes transform ${T}(z)$ of the limiting eigenvalue distribution
of $Q_{\ul{k}}(\xi_n)/n$ with i.i.d.\ matrices $\xi_n\in E_{\ul{k}}$ is given for $z\in\C^+$ as
\[
{T}(z)=-\frac{1}{v}-\frac{c}{zW_{\kappa,\gamma}(-\frac{cv}{z})}-\frac{c\gamma+1-c}{z}
=\frac{\exp_\kappa\bigl(W_{\kappa,\gamma}(-vc/z)\bigr)-1}{v}-\frac{1-c}{z}.
\]
\end{Corollary}
The properties of absolute continuity  and support of the limiting measure can be derived  analogously  to
those obtained in Theorem~\ref{th:WishartProfile}
for $c=1$.
\begin{proof}
%\sout{We shall use the same method as before so that}
We construct a variance profile $\sigma$ from $E_{\ul{k}}$ likely to~\eqref{def:varianec profile}.
We embed the rectangular matrix $\xi_n\in E_{\ul{k}}$ in a square matrix
$Y(\xi_n)=\pmat{0&\xi_n\\ {}^{t}\xi_n&0}$, and set $V_n=\set{Y(\xi_n)}{\xi_n\in E_{\ul{k}}}$. 
Set $\ds p'=\lim_{n\to\infty}\frac{n}{n+N_n}=\frac{1}{1+m_1+m_2(1-c)}$.
Let $\sigma$ be a function $[0,1]\times[0,1] \to\R_{\ge0}$ defined by
\[
\sigma(x,y)=\begin{cases}
v&(\text{$x<cp'$ and $y\ge p'+m_1x$})\\
v&(\text{$x\ge p'$ and $0\le y\le \min\{(x-p')/m_1,cp'\}$}),\\
0&(\text{otherwise}).
\end{cases}
\]
Then, we can show that $\sigma$ is the variance profile of $V_n$.
On the other hand,
let us consider a subspace $E'_{\ul{k}}:=\set{\xi=\pmat{\eta\\ \zeta}\in E_{\ul{k}}}{\zeta=0}$ of $E_{\ul{k}}$,
and let $V'_n=\set{Y(\xi_n)}{\xi_n\in E'_{\ul{k}}}$.
Then, 
$\sigma$ is also the variance profile of $V'_n$.
Thus, we consider equivalently the limiting eigenvalue distribution of $V'_n$,
and  that of covariance matrices on $E'_{\ul{k}}$.
If $\xi_n=\pmat{\eta_n\\0}\in E'_{\ul{k}}$, then 
$Q_{\ul{k}}(\xi_n)=\pmat{\eta_n{}^{\,t\!}\eta_n&0\\0&0}$,
and thus it is enough to study the limiting eigenvalue distribution of $\eta_n{}^{\,t\!}\eta_n$.
%\sout{ Since  it}
The variance profile of $\eta_n{}^{\,t\!}\eta_n$ 
has a trapezoidal form \eqref{def:varianec profile} 
(illustrated by \eqref{eq:fig sigma})  
with parameters $\alpha=m_1$ and $p=\lim_n\frac{a_n}{a_n+N_n}=\frac{c}{c+m_1+m_2(1-c)}$.
Applying Proposition~\ref{prop:Stieltjes of quad Sym},
we see that
the corresponding Stieltjes transform $T_1(z)$ is given by
\[
T_1(z)=T_{\kappa,\gamma}(z)\quad\text{with}\quad
\kappa=\frac{1}{1-m_1},\quad
\gamma=\frac{2p-1}{p}=\frac{c-m_1-m_2(1-c)}{c}.
\]
In general, for two symmetric matrices $A_i$ $(i=1,2)$ of size $n_i$,
the Stieltjes transform $S(z)$ of $\mathrm{diag}(A_1,\,A_2)/(n_1+n_2)$ can be described by using the Stieltjes transforms
$S_i(z)$ of $A_i/n_i$ $(i=1,2)$ as
\[
S(z)=S_1\left(\frac{n_1+n_2}{n_1}z\right)+S_2\left(\frac{n_1+n_2}{n_2}z\right)\quad(z\in\C^+).
\]
In our situation, we have $(n_1,n_2)=(a_n,b_n)$ and 
$(A_1,A_2)=(\eta_n{}^{\,t\!}\eta_n,0)$.
Hence, we have $S_2(z)=-\frac{1}{z}$
and
$S_1(z)$ is the Stieltjes transform of $\eta_n{}^{\,t\!}\eta_n/{a_n}$
so that $\lim_{n\to\infty} S_1(z)=T_1(z)$.
Thus, taking the limit $n\to\infty$, 
we see that
the limiting 
%\sout{covariance} 
Stieltjes transform $T(z)$ {corresponding to} $E'_{\ul{k}}$, and hence {to} $E_{\ul{k}}$ is given as
\[
\begin{array}{r@{\ }c@{\ }l}
T(z)
&=&
\ds
T_1\left(\frac{z}{c}\right)+S_2\left(\frac{z}{1-c}\right)
=
T_{\kappa,\gamma}\left(\frac{z}{c}\right)-\frac{1-c}{z}\\[0.8em]
&=&
\ds
-\frac{1}{v}-\frac{c}{zW_{\kappa,\gamma}(-vc/z)}-\frac{c\gamma+1-c}{z},
\end{array}
\]
whence we obtain the corollary.
\qed
\end{proof}

\begin{figure}
\begin{tabular}{ccc}
    %\begin{minipage}{100pt}
    %   \centering
    %   \includegraphics[scale=0.2]{Graphs/Wishart(a0n4000p2).jpg}
    %%   \caption{Simulation for $\alpha=0$}
    %   \label{fig:Wishart1}
    %\end{minipage}
    %&
    \begin{minipage}{0.3\textwidth}
        \centering
        \includegraphics[scale=0.25]{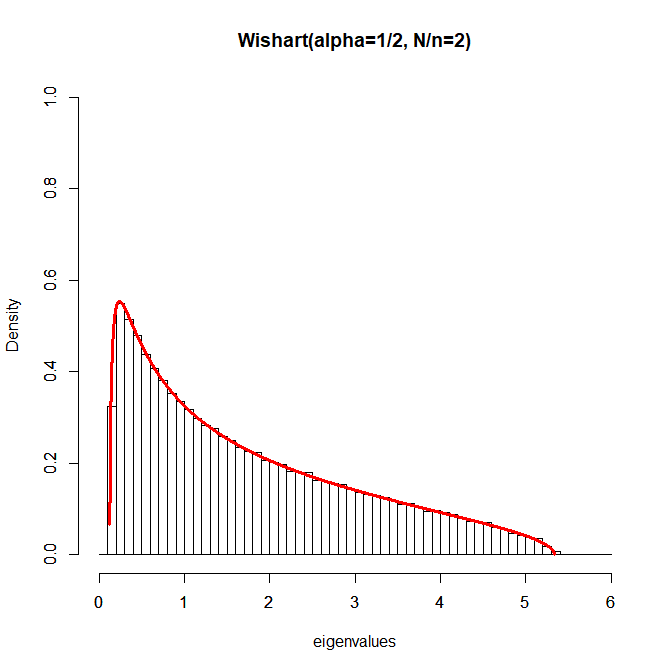}
        \caption{Simulation for $\alpha=1/2$}
        \label{Wishart2}
    \end{minipage}
    &
    \begin{minipage}{0.3\textwidth}
        \centering
        \includegraphics[scale=0.25]{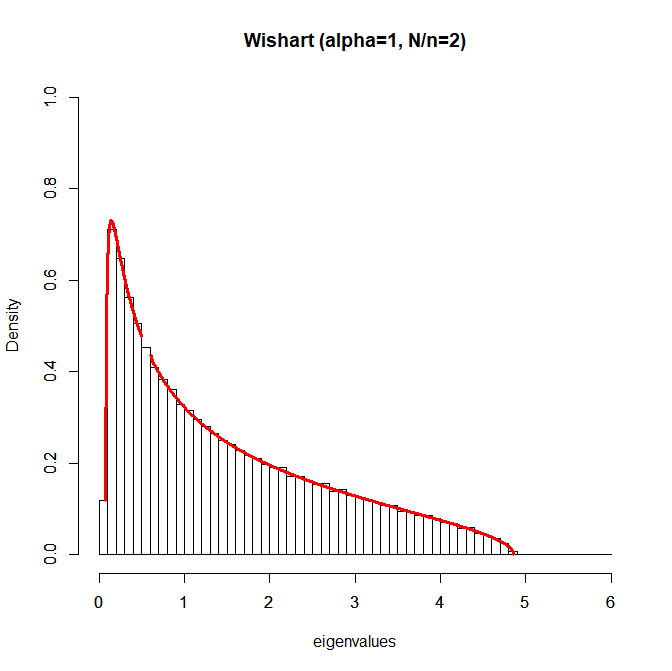}
        \caption{Simulation for $\alpha=1$}
        \label{fig:Wishart3}
    \end{minipage}
    &
    \begin{minipage}{0.3\textwidth}
        \centering
        \includegraphics[scale=0.25]{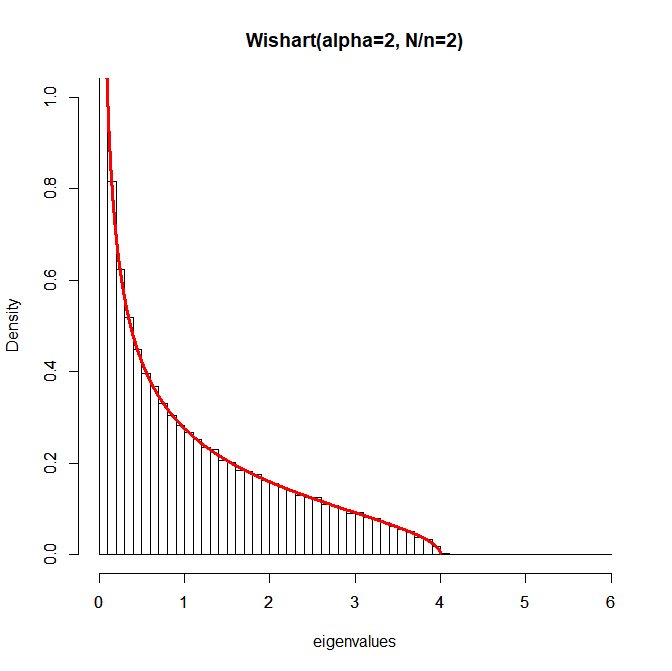}
        \caption{Simulation for $\alpha=2$}
        \label{Wishart4}
    \end{minipage}
\end{tabular}
\end{figure}

\begin{Remark}
\label{remark}
{\rm
In the Figures \ref{Wishart2}-\ref{Wishart4}
we present simulations of $\ul{k}$-indexed
Wish-art ensembles  $X_n=Q_{\ul{k}}(\xi_n)$
on the symmetric cone $\mathrm{Sym}(n,\R)^+$
(i.e. $c=1$), 
for $n=4000$ and $N=|\ul{k}|=2n$ with parameters $\alpha=m_1=1/2$, $1$ and $2$, respectively. 
We have $\gamma=-1$ and $\kappa=2,\infty, -1$ respectively.
The red line is the graph of $d(x)$ generated by the R program from its Stieltjes transform given 
in Corollary ~\ref{prop:support}.
In two first cases, the limiting density $d(x)$ is continuous on $\R$ with compact support contained in $(0,\infty)$.
The last case $(\kappa,\gamma)=(-1,-1)$ 
corresponds to $(\kappa',\gamma')=(1,0)$ 
which is the classical Wishart case with $C=1$.
Thus its density explodes to $\infty$ at 0.
}
\end{Remark}
\out{\bl{
The following histograms are simulations of Wishart ensembles for $n=4000$ and $b_n=1$ with parameters $m_1=0,\frac12,1,2$ and $|\ul{k}(n)|=2n$.
\sout{These are samples.}
}}

%%%%%%%%%%%%%%%%%%%%%%%%%%%%
\begin{Remark}\label{otherWishart}
{\rm
Let $Y_n$ be a rectangular $n\times p$
i.i.d.\ matrix with variance profile $\sigma(x,y)$,
and assume that $\lim_{n\to\infty}p/n = c$.
In papers \citet{HLN,HLN-AHP,HLN-AAP515}
a functional equation 
$\tau(u,z)=
\bigl(-z +\int_0^1 \sigma(u,v)
\bigl(1+c\int_0^1 \sigma(x,v)\tau(x,z)dx\bigr)^{-1}dv\bigr)^{-1}
$
is given to get  the limiting Stieltjes  transform
$f(z)$   for the   rescaled random matrices $Y_nY_n^*$, 
as the
integral $\int_0^1 \tau(u,z)du$. 
This equation appears in 
\citet{GirkoBook}  in the setting of Gram matrices based on Gaussian
fields
(cf.\ \citet[Remark 3.1]{HLN-AHP}).

However, thanks to symmetry, solving the equations ~\eqref{def:difEqQW} resulting from Theorem~\ref{theo:variance profile method} is easier than solving the last functional-integral equation for $\tau(u,z)$.
Therefore
we opted for variance profile method for Gaussian and Wigner ensembles as the main tool  of studying Wishart ensembles of Vinberg matrices.
}
\end{Remark}

\section{Wigner and Wishart Ensembles related to Generalized Vinberg cones}
\label{sect:5}

%(SLIGHTLY REWRITTEN BY PIOTR)

In this section,
we consider the dual cone $\Vin$ of $\dVin$, 
which is realized as a minimal matrix form in the sense of \citet{YN2015} as follows.
Let $\mathcal{V}_n$ be a subspace of $\mathrm{Sym}(a_n(b_n+1),\,\R)$ defined by
\begin{equation}
\label{eq:elementVn}
\mathcal{V}_n:=
\set{\mathrm{diag}\left(\pmat{x&y_1\\{}^{t\!}y_1&d_1},\dots,\pmat{x&y_{b_n}\\{}^{t\!}y_{b_n}&d_{b_n}}\right)}{
\begin{array}{l}
x\in\mathrm{Sym}(a_n,\R),\\ 
y_1,\dots,y_{b_n}\in\R^{a_n},\\
d_1,\dots,d_{b_n}\in\R
\end{array}
}.
\end{equation}
Then, the dual cone $\Vin$ is described as $\Vin:=\mathcal{V}_n\cap\mathrm{Sym}(a_n(b_n+1),\,\R)^+$.

We consider 
Wigner Ensembles $V_n\in\mathcal{V}_n$ 
and quadratic Wishart Ensembles $X_n\in\Vin$
as those in the sense of $\mathrm{Sym}(a_n(b_n+1),\,\R)$.
Assume that $\lim_{n\to+\infty}a_n=\infty$.
By the theory of lower rank perturbation 
(see \citet[{\S}2.4.1]{Tao},  for example), 
the study of 
eigenvalue distributions
of these ensembles
boils down to the study of the eigenvalue distributions of $x$ and, after suitable normalization,
the limiting eigenvalue distributions of $V_n$ and $X_n$ are the same as for $x\in \mathrm{Sym}(a_n,\,\R)$.

This essential difference in the
Random Matrix Theory for the cones
$\Vin$ and $\dVin$ may be explained by
a substantial difference between the cones
$\Vin$ and $\dVin$ in terms of
numbers of \textit{sources} in the sense of~\citet{YN2015}.
In the case $\dVin$, there is only one source so that $\dVin$ can be realized in a usual matrix form.
On the other hand, 
$\Vin$ has $b_n$ sources so that 
%\sout{we need} 
$b_n$ copies of a usual matrix form appear.

\section{Acknowledgements}
%(Professor H.\ Ishi, Professor J. Najim, Grant of H. Nakashima)
The authors would like to express their sincere gratitude to Professor H.\ Ishi for
strong encouragement and invaluable comments on this work.
The authors are very  grateful to Professor J.\ Najim for significant
methodological and bibliographical indications for this work 
and to Professor C.\ Bordenave for discussions on Theorem
\ref{theo:variance profile method}.
The authors thank the Scientific Committee of the CIRM Luminy 2020
Conference "Mathematical Methods of Modern Statistics" (MMMS 2)  for giving the possibility to present this work. We thank numerous participants  of MMMS 2
for their comments and remarks.

This research was carried out while the first author spent the winter semester of 2018 at 
Laboratoire de Math{\'e}matiques LAREMA
under the support of Grant-in-Aid for JSPS fellows (2018J00379).

%%%%%%%%%%%%%%%%%%%%%%%%%%%%%%%%%%%%%%%%%%%%%%%%%%%%%%%%%%%%

%%%%%%%%%%%%%%%%%%%%%%%%%%%%%%%%%%%%%%%%%%%%%%%%%%%%%%%%%%%%%

%%%%%%%%%%%%%%%%%%%%%%%%%%%%%%%%%%%%%%%%%%%%%%%%%%%%%%%%%%%%%%%%%%%%%%%%%%%%%%%%%%%%%%%%%%%%%%%%%%
\out{

}
%%%%%%%%%%%%%%%%%%%%%%%%%%%%%%%%%

\end{document}